\documentclass [12pt,letterpaper]{article}

\usepackage{amsmath}		% Extra math definitions
\usepackage{graphics} 		% PostScript figures
\usepackage{setspace}		% 1.5 spacing
\usepackage{longtable}          % Tables spanning pages
\usepackage{cleveref}
\usepackage{verbatim}
\usepackage{bbm}
\usepackage{fancyhdr}      % (Ali) To change the heading in the appendix

\usepackage{todonotes}

\usepackage{color}
\usepackage{pdfpages}
\usepackage{graphicx}
\usepackage{epsfig} % for postscript graphics files
\usepackage{times} % assumes new font selection scheme installed
\usepackage{amsmath} % assumes amsmath package installed
\usepackage{amsthm}
\usepackage{amssymb}  % assumes amsmath package installed
\usepackage{amsfonts}
\usepackage{dsfont}
\usepackage{tikz,tabstackengine,amsmath}
\usetikzlibrary{shapes,shadows,arrows}
\usetikzlibrary{shapes,shapes.geometric,arrows.meta,fit,calc,positioning,automata}
\usepackage{pgfplots}
\usepgfplotslibrary{groupplots}
\usepackage{subfig}
\usepackage{mathtools}
\usepackage{adjustbox}

\usepackage [english]{babel}
\usepackage [autostyle, english = american]{csquotes}
\MakeOuterQuote{"}

\usepackage{bm}

\usepackage[utf8]{inputenc}

\usepackage{enumitem}

\usepackage[noadjust]{cite}

\newtheorem{theorem}{Theorem}
\newtheorem{lemma}[theorem]{Lemma}

\theoremstyle{definition}

\theoremstyle{remark}
\newtheorem{remark}{Remark}

\theoremstyle{assumption}
\newtheorem{assumption}{Assumption}

\newcommand{\mb}{\mathbb}
\newcommand{\mc}{\mathcal}

\newcommand{\mdsR}{\mathds{R}}
\newcommand{\mfM}{{\mathfrak{M}}}
\newcommand{\mfN}{{\mathfrak{N}}}
\newcommand{\mfK}{{\mathfrak{K}}}
\newcommand{\mfNo}{{\mfN_{0}}}

\begin {document}

%===== Title page
\title{Mean Field Game Systems with Common Noise and Markovian Latent Processes}
\author{Dena Firoozi, Peter E. Caines, Sebastian Jaimungal}
\date{}

\maketitle

%\label{chapter:MFG-Latent} 
\section{Abstract}
In many stochastic games stemming from financial models, the environment evolves with latent factors and there may be common noise across agents' states. Two classic examples are: (i) multi-agent trading on electronic exchanges, and (ii) systemic risk induced through inter-bank lending/borrowing. Moreover, agents' actions often affect the environment, and some agent's may be small while others large. Hence sub-population of agents may act as minor agents, while another class may act as major agents. To capture the essence of such problems, here, we introduce a general class of non-cooperative heterogeneous stochastic games with one major agent and a large population of minor agents where agents interact with an observed common process impacted by the mean field. A latent Markov chain and a latent Wiener process (common noise) modulate the common process, and agents cannot observe them. We use filtering techniques coupled with a convex analysis approach to (i) solve the mean field game limit of the problem, (ii) demonstrate that the best response strategies generate an $\epsilon$-Nash equilibrium for finite populations, and (iii) obtain explicit characterisations of the best response strategies.

\section{Introduction}
%\seb{add references to Huyen Pham, Rene Aid for commodities Ulrich Horst and Charles Lehalle and Nourin on MFG in Algo}
 Mean Field Game (MFG) systems theory establishes the existence of approximate Nash equilibria and the corresponding individual strategies for stochastic dynamical systems in games involving a large number of agents \cite{CainesEncyc2014, HuangCIS2006, HuangCDC2006}. The equilibria are termed $\epsilon$-Nash equilibria and are generated by the local, limited information feedback control actions of each agent in the population. The feedback control actions constitute the best response of each agent with respect to the behaviour of the mass of agents. Moreover, the approximation error, induced by using the MFG solution, converges to zero as the population size tends to infinity.                       \begin{comment}                         MFG theory offers a framework in which each agent interacts with the aggregate effect of all other agents and optimizes their individual cost functional based on (i) local information of its own state, and (ii) information on the overall population state, i.e. the mean field.
\end{comment}
The analysis of this set of problems originated in \cite{HuangCDC2003a, HuangCDC2003b, HuangCIS2006, HuangTAC2007} (see \cite{CHM17}, \cite{CainesEncyc2014}), and independently in \cite{Lasry2006a, Lasry2006b, Lasry2007}. In \cite{Huang2010} and \cite{Huang2012} the authors analyse and solve the linear quadratic systems case where there is a major agent (whose impact does not vanish in the limit of infinite population size)  together with a population of minor  agents (where each agent  has individually asymptotically negligible effect). Major-minor problems lead to stochastic mean fields; and by extending the state space for such systems, \cite{Huang2012} and \cite{NourianSiam2013} establish the existence of $\epsilon$-Nash equilibria together with the individual agents' control laws that yield the equilibria. Following this work, the partially observed MFG theory for nonlinear and linear quadratic systems with major and minor agents is developed in \cite{Kizilkale2013, Kizilkale2014, KizilkaleTAC2016, SenCDC2014, SenACC2015, SenSiam2016, FirooziCDC2015} where the agents' states are partially obfuscated by Brownian noise. 

The MFG approach has been applied recently in a number of  financial modeling problems. Here we highlight a few important contributions. An optimal execution problem in algorithmic trading with the linear models in \cite{JaimungalBook2015} was formulated as for the nonlinear major-minor (MM) MFG model in \cite{Mojtaba2015}. The completely observed and partially observed major minor linear quadratic MFG theory was first applied to an optimal execution problem with linear models (see e.g., \cite{JaimungalBook2015}) in \cite{FirooziCDC2016}, and subsequently in \cite{FirooziCDC2017,FirooziIFAC2017, FirooziISDG2017}. Optimal stopping and switching time problems are addressed for competitive market participants in \cite{FirooziPakniyatCainesCDC2017}. A mean field game of controls is applied to optimal trading problems in \cite{Cardaliaguet2018,Lehalle2019}. Financial markets driven by latent factors are modeled in the MFG framework in \cite{JaimungalPhil2018, casgrain2018mean}. The problem of optimal portfolio liquidation using MFG formulation was addressed in \cite{Horst2018}.  \cite{Aid2017} investigates investment decisions in the electricity market using   linear-quadratic  McKean-Vlasov control problems with random coefficients.

\textit{Problem description:} In this paper, an MFG framework is considered where there exist one major agent and a large number of minor agents which are subject to linear dynamics and quadratic cost functionals. Each agent interacts with other agents in the system through the coupling in their cost functional with a common process. The common process is modulated by a latent (unobservable) Markov  process, the major agent's state, the major agent's control action, the average state of all minor agents, the average  control action of all minor agents, and a common Wiener process. The objective is to characterize the best response strategies which yield an $\epsilon$-Nash equilibrium.  

\textit{Motivation:} 
Financial and economic systems (among others) are often driven by latent factors, and these latent factors also affect the cost (profit) functional of the traders involved. Moreover, the agents in these system are often acting in a non-cooperative manner, and hence playing a large stochastic game with one another; while they may control aspects of the system, they are also at the whim of factors they cannot control or observe. For example, in optimal execution problems (where traders aim to sell or buy shares of an asset), all traders are subject to the same asset price process and must make their trading decisions based on the observed price. The asset price dynamics may be driven by a common Wiener process, which accounts for so-called noise traders who are considered non-sophisticated and uninformed about future price movements. In addition, the effect of unobserved latent factors, other than the major agent's trading action and the aggregate impact of minor agents' trading actions, are important to incorporate (see e.g. \cite{JaimungalPhil2016}, \cite{JaimungalPhil2018} ) in specifying the best response trading strategies and resulting $\epsilon$-Nash equilibrium. 

% The asset price  may be modeled by each agent as an SDE involving with a common Wiener process. Latent processes may also be used to model the effect of unobserved factors on the price dynamics, other than the major agent's and the aggregate impact of minor agents' trading actions (see \cite{JaimungalPhil2016}).

\textit{Methodology:} Although latent processes are not directly observable, the information provided from the realized trajectories of the common process and the evolution of system's aggregate state (mean field) can be used to obtain posteriori estimates, and to subsequently partially predict future behavior of the common process \cite{JaimungalPhil2016}. Certain versions of such problems can then be recast as MFG systems with a common noise. A variation of this type of MFG system has been investigated in \cite{CarmonaDelarue2016}, where the case of correlated randomness in a nonlinear setting is analyzed. Here we utilize a different approach in order to address the existence of a latent process together with the common noise.
Specifically,  we treat the common process as a major agent and further extend the Major-Minor LQG MFG analysis of \cite{Huang2010} to incorporate a latent process in the dynamics. Moreover, we develop a novel convex analysis approach to obtain the best response strategies for all agents that yield  an $\epsilon$-Nash equilibrium. We presented the initial results of this work in \cite{FJCCDC2018}.

The rest of the paper is organized as follows. Section \ref{sec:setup} introduces a class of major-minor MFG problems with  a common process as well as a latent process. The MFG formulation of the problem is then presented in Section \ref{sec:MFG formulation}. Concluding remarks are made in Section \ref{sec:conclusion}.

\section{Major-Minor Mean Field Game Systems with a Common Process} \label{sec:setup}

\subsection{Dynamics: Finite Population}
We consider a large population of $N$ minor agents and a  major agent, where the agents are coupled through their individual cost functionals with a common process. 
\subsubsection{Major and Minor agents}
The underlying dynamics of the major and minor agents are assumed to be given, respectively, by
\begin{align}
dx_t^0 &= [A_0x_t^0 + B_0 u_t^0 + b_0(t)]dt + \sigma_0 dw^0_t, \label{GSCMFGMajorIncomDyn}\\
dx_t^i &= [A_k x_t^i + B_k u_t^i + b_k(t)]dt + \sigma_k dw^i_t, \label{GSCMFGMinorIncomDyn} 
\end{align}
where $t \in [0,T]$, $ i \in\mfN:=\{1,\dots,N\}$, $N< \infty$, and the subscript $k \in \mfK:=\{1,\dots,K\},\, K\leq N$, denotes the type of a minor agent. Here, $x^i_t \in \mdsR^n$, $i\in\mfN_0:=\mfN\cup\{0\}$, are the states, $u^i_t \in \mdsR^m$, $i\in\mfNo$ are the control inputs, $\lbrace w^i_t,~ i\in\mfNo \rbrace$ denotes $(N+1)$ independent standard Wiener processes in $\mdsR^r$, where $w_i$ is progressively measurable with respect to the filtration $\mathcal{F}^w \coloneqq (\mathcal{F}_t^{w})_{t\in [0,T]}$. All matrices in \eqref{GSCMFGMajorIncomDyn} and \eqref{GSCMFGMinorIncomDyn} are constant and of appropriate dimension; vectors $b_0(t)$, and $b_k(t)$ are deterministic functions of time.

\begin{assumption} \label{ass:IntialStateAss} 
The initial states $\lbrace x^i_0,~ i \in \mathfrak{N}_0\rbrace$ are identically distributed and mutually independent and also independent of $\mathcal{F}^{w}$; $\mathbb{E}[ w^i_t (w^i_t)^T] = \Sigma$. Moreover, $\mathbb{E}x^i_0= 0$, and $\mathbb{E}\Vert x^i_0\Vert^2 \leq C < \infty $, with $\Sigma$ and $C$ independent of $N$, for all $i\in\mfNo$.
\end{assumption}

\subsubsection*{Minor Agents Types:}

Minor agents are given in $K$ distinct types with $1 \leq K < \infty$. The notation 
\begin{align*}
\Psi_k \triangleq \Psi(\theta_i), \quad \theta_i = k 
\end{align*}
is introduced where $\theta_i \in \Theta$, with $\Theta$ being the parameter set, and $\Psi$ may be any dynamical parameter in \eqref{GSCMFGMinorIncomDyn} or weight matrix in the cost functional  \eqref{GSCMFGlatent-minorcost}. The symbol $\mathcal{I}_k$ denotes
\begin{align}
\mathcal{I}_k = \lbrace i : \theta_i = k,~ i \in \mathfrak{N} \rbrace , \quad k \in \mfK\nonumber
\end{align}
where the cardinality of $\mathcal{I}_k$ is denoted by $N_k = |\mathcal{I}_k|$. Then, $\pi^{N} = (\pi_{1}^{N},...,\pi_{K}^N),~ \pi_k^N = \tfrac{N_k}{N}, ~ k \in \mfK$, denotes the empirical distribution of the parameters $(\theta_1,...,\theta_N)$ sampled independently of the initial conditions and Wiener processes of the agents $\mathcal{A}_i, i \in \mathfrak{N}$. The first assumption is as follows.
\begin{assumption} \label{EmpiricalDistAss}
There exists $\pi$ such that $\mbox{lim}_{N \rightarrow \infty} \pi^N = \pi $ a.s. 
\end{assumption}
\subsubsection{Common Process: Finite Population}
We consider the systems where the major agent and any minor agent $\mc{A}_i,\, i \in \mathfrak{N}$, observe a common stochastic process $y_t$, where both the state and common process $y_t$ appear in an agent's cost functional as will be introduced in Section \ref{sec:costs}. The common process $y_t \in \mdsR^n$ is governed by
%\begin{subequations}\label{grp}
\begin{align} \label{GSCMFGminorComDynNonEst}
&dy_t = dy_t^L + (F u_t^{(N)} dt + F_0 u^0_t + H x^{(N)}_t + H_0 x^0_t)dt, 
\end{align}
%\end{subequations}
where $y_t^L \in \mdsR^d$ evolves as in
\begin{align}
&dy_t^L  = f(t,y_t^L, \Gamma_t) dt + \sigma dw_t. \label{GSCMFGminorComDynNonEst-unimpacted}
\end{align}
In  \eqref{GSCMFGminorComDynNonEst-unimpacted}, the process $\Gamma \coloneqq (\Gamma_t)_{t \in [0,T]}$ denotes a latent continuous Markov chain process with $\Gamma_t \in \{ \gamma_j,\,  j \in \mathfrak{M}\},\, \mathfrak{M}=\{1,\dots,M\},\, M< \infty$; the vector $f(t,y_t^L, \Gamma_t)$ denotes a deterministic nonlinear function of $t$, $y^L$, and $\Gamma$; $w_t \in \mdsR^r$ denotes a latent Wiener process independent of $\{ w^i_t,\, i \in \mathfrak{N}_0 \}$, and the matrices $F$, $F_0$, $H$, $H_0$, and $\sigma$ are deterministic, constant and of appropriate dimensions. Moreover, by substituting \eqref{GSCMFGminorComDynNonEst-unimpacted} in \eqref{GSCMFGminorComDynNonEst}, it is evident that the common process $y_t$ is impacted by
\begin{flushleft}
\begin{enumerate}
\item[1)] a latent Markov chain process $\Gamma_t$,% not directly observable to each agent,
\item[2)]  the major agent's state $x^0_t$,
\item[3)] the major agent's control action $u^0_t$, 
\item[4)] the average state of minor agents, i.e. $x^{(N)}_t = \frac{1}{N}\sum_{i=1}^N x^i_t$,
\item[5)] the average control action of minor agents, i.e. $u^{(N)}_t = \frac{1}{N}\sum_{i=1}^{N} u^i_t$,
\item[6)] a Wiener process $w_t \in \mdsR^r$ independent of   $w^i_t,~ i\in\mfNo$.
 \end{enumerate}
 \end{flushleft}
%\begin{assumption}\label{ass:CommonObs} The major agent $\mc{A}_0$ and each minor agent $\mc{A}_i,\, i \in \mathfrak{N} < \infty$, have no observations on the continuous Markov chain process $\Gamma \coloneqq (\Gamma_t)_{t \in [0,T]}$. 
%\end{assumption}
\begin{assumption} \label{ass:majorObs}
The major agent $\mc{A}_0$ completely observes its own state $x^0_t$, and the common process $y_t$.
\end{assumption}

\begin{assumption} \label{ass:minorObs}
Each minor agent $\mc{A}_i,\; i \in \mathfrak{N}$, completely observes its own state $x^i_t$, the major agent's state $x^0_t$, and the common process $y_t$. 
\end{assumption}
We again emphasize that the latent processes $\Gamma_t$ and $w_t$ are not directly observed by the agents $\mc{A}_i, \, i \in \mathfrak{N}_0$. However, each agent may obtain their posteriori estimates based on its complete observations on the common process $y_t$. We refer to the latent Wiener process as the common noise process in this work.

\subsubsection{Control $\sigma$-Fields}

For any finite $T$, We denote (i) the natural filtration generated by the $i$-th minor agent's state $(x^i_{t})_{t\in[0,T]}$ by $\mathcal{F}^i\coloneqq (\mathcal{F}_{t}^i)_{t\in [0,T]}$, $ i\in\mfN$, (ii) the natural filtration generated by the major agent's state $(x^0_{t})_{t\in[0,T]}$ by $\mathcal{F}^0:=(\mathcal{F}^0_t)_{t\in[0,T]}$ , (iii) the natural filtration generated by the states of all agents $((x_t^i)_{i\in\mfNo})_{t\in[0,T]}$ by $\mathcal{F}:=(\mathcal{F}_{t})_{t\in[0,T]}$, (iv) the natural filtration generated by  $(\Gamma_{t}, w_{t})_{t\in[0,T]}$ by $\mc{G}:=(\mc{G}_t)_{t\in[0,T]}$ , and (v) the natural filtration generated by $(y_{t})_{t\in[0,T]}$ by $\mc{F}^y \coloneqq (\mc{F}_t^y)_{t\in[0,T]}$. %Finally $\mc{F}_0$ is defined to be the increasing family of $\sigma$-fields generated by the initial states $(x^i_0; ~i\in\mfN)$.

%Next, we introduce three admissible control sets. First,  $\mc{U}^{0,L}$ denotes the set of controls based on the local information set of the major agent and consists of the set of linear feedback control laws $u^0$ that are adapted to  $\mc{F}^0$ such that $\mb{E}[\int_0^T u_t^{0\intercal}u_t\, dt] < \infty$. Second, for each $i \in \mfN$, $\mc{U}^{i,L}$ denotes the set of controls based on the local information set of  $\mc{A}_i$ and consists of the linear feedback control laws adapted to the filtration $\mc{F}^{i,r} \coloneqq (\mc{F}_{t}^{i,r})_{t\in\mfT}$, where $\mc{F}^{i,r} \coloneqq \mc{F}^i \vee \mc{F}^0$, such that  $\mb{E}[\int_0^T u_t^{i\intercal}  u_t^i\, dt] < \infty$. Third,  $\mc{U}_g^{N,L}$ denotes the set of linear feedback controls $u$ that are adapted to the general filtration $\mc{F}^g:=(\mc{F}^g _{t})_{t\in\mfT}$,  $\mc{F}^g:= \lor_{i \in \mfN_0} \mc{F}^i$, such that $\mb{E}[\int_0^T u_t^\intercal  u_t\,dt]< \infty$.

Next, we introduce three admissible control sets. First, $\mathcal{U}^0$ denotes the set of feedback control laws adapted to the filtration $\mc{F}^{0,r}\coloneqq (\mc{F}_t^{0,r})_{t \in [0,T]}$, where $\mc{F}^{0,r} \coloneqq \mathcal{F}^0 \vee \mathcal{F}^y$, such that $\mb{E}[\int_0^T u_t^{0\intercal}u_t\, dt] < \infty$. Second, $\mathcal{U}^i$, $i\in\mfN$, denotes the set of control laws adapted to the filtration $\mathcal{F}^{i,r} \coloneqq (\mathcal{F}_{t}^{i,r})_{t\in[0,T]}$, where $\mc{F}^{i,r} \coloneqq \mathcal{F}^i \vee \mathcal{F}^0 \vee \mathcal{F}^y$, $i\in\mfN$, such that  $\mb{E}[\int_0^T u_t^{i\intercal}  u_t^i\, dt] < \infty$. Third,  $\mathcal{U}_g^{N}$ is adapted to the general filtration $\mathcal{F}^{g}:=(\mathcal{F}_{t}^{g})_{t\in[0,T]}$, where $\mc{F}^g \coloneqq \mathcal{F} \vee \mathcal{F}^y \vee \mathcal{G}$, such that $\mb{E}[\int_0^T u_t^\intercal  u_t\,dt]< \infty$. %We note in passing the significant differences between the information structures specified here and those in the cooperative game theory literature \cite{NayyarTAC2014}.
\begin{assumption}
[Major Agent's Linear Control Laws] For major agent $\mc{A}_0,$ the set of control laws $\mc{U}^{0,L} \in \mc{U}^0,$ is
defined to be the collection of linear feedback control laws adapted to $\mc{F}^{0,r}$.
\end{assumption}
\begin{assumption}
[Minor Agent's Linear Control Laws] For each
minor agent $\mc{A}_i,\, i \in \mathfrak{N}$, the set of control laws $\mc{U}^{i,L} \in \mc{U}^i, \, i \in \mathfrak{N}$, is
defined to be the collection of linear feedback control laws adapted to $\mc{F}^{i,r},\,i \in \mathfrak{N}$. 
\end{assumption}

\subsection{Cost Functionals: Finite Population } \label{sec:costs}
Given the vector $z^0_t$ as 
\begin{align*}
z^0_t = \left[ 
\begin{array}{c}
y_t^\intercal, \;
x_t^{0\intercal}
\end{array} \right]^\intercal.
\end{align*}
the major agent's cost functional to be minimized is formulated by 
\begin{equation}\label{GSCMFGlatent-majorcost}
J_0(u^0, u^{-0}) = \tfrac{1}{2}\mb{E} \bigg [  z_T^{0\intercal} G_0 z_T^0 + \int_0^T \Big \{z_s^{0\intercal} Q_0 z_s^0 + 2\,z_s^{0\intercal} N_0 u_s^0 + u_s^{0\intercal} R_0 u^0_s \Big \} ds \bigg],
\end{equation}
where $u^{-0}=(u^1,u^2,...,u^{N})$. 

\begin{assumption} \label{ConvexityCondMajor}
$G_0$, $Q_0$, and $R_0$ are symmetric matrices of appropriate dimension, and we have $G_0 \geq 0$, $Q_0 - N_0 R_0^{-1} N_0^\intercal \geq 0$, and  $R_0 >\delta I$ for some $\delta>0$ (see \cite{XYZBook1999}).  
 
\end{assumption}
Similarly, to define the minor agents' cost functional, let the vector $z^i_t, \, i\in\mfN$, denote\footnote{For the minor agent, in $z_t^i$ we order $y_t$ second, while for the major agent in $z^0_t$ we order $y_t$  first. This is done for convenience when we extend the minor agents' state  in Section \ref{sec:minorInf}.}
\begin{align*}
z^i_t = \left[ \begin{array}{c}
x_t^{i\intercal},\;
y_t^\intercal
\end{array} \right]^\intercal. 
\end{align*}
Then, minor agent $\mathcal{A}_i,\, i \in \mfN$, of type $k \in \mfK$,  cost functional to be minimized is

\begin{equation}\label{GSCMFGlatent-minorcost}
%\begin{split}
J_i(u^i, u^{-i}) = \tfrac{1}{2}\mb{E} \bigg [  z_T^{i\intercal} G_k z_T^i + \int_0^T \Big \{z_s^{i\intercal} Q_k z_s^i + 2 z_s^{i\intercal} N_k u_s^i + u_s^{i\intercal} R_k u^i_s \Big \} ds \bigg], 
%\end{split}
\end{equation}
$k\in\mfK$, where $u^{-i} = (u^0,...,u^{i-1}, u^{i+1},...,u^{N})$.%\seb{you need to identify agent-$i$ as being in sub-population-$k$}
\begin{assumption}\label{convexityCondMinor}
 $G_k$, $Q_k$, and $R_k$ are symmetric matrices of appropriate dimension, and we have $G_k \geq 0$, $Q_k - N_k R_k^{-1} N_k^\intercal \geq 0$, and  $R_k >\delta I$ for some $\delta>0$ and $\forall k \in \mfK$ (see \cite{XYZBook1999}).  
\end{assumption}

%%%%%%%%%%%%%%%%%%%%%%%%%%%%%%%%%%%%%%%%%%%%%%%%%%%%%%%%%%%%%%%%%%%%%%%%%%%%%%%%
\section{Major-Minor LQG Mean Field Games Approach}\label{sec:MFG formulation}
In the mean field game methodology with a major agent \cite{NourianSiam2013}, \cite{Huang2010}, the problem is first solved in the infinite population case where the average terms in the finite population dynamics and cost functional of each agent are replaced with their infinite population limit, i.e. the mean field. For this purpose, the major agent's state is extended with the mean field, while the minor agent's state is extended with the major agent's state, and the mean field; this yields stochastic optimal control problems for each agent linked only through the major agent's state and mean field. Finally the infinite population best response strategies are applied to the finite population system which yields an $\epsilon$-Nash equilibrium. 
%\seb{shouldn't we  hyphenate ``major-minor'' throughout?}\dena{I did the replacement} 

To address major-minor mean field game systems involving a common process which, in turn, is impacted by a latent Markov chain process, the following steps are followed. We first note that the common process in this work represents an extended form of common noise in \cite{CarmonaDelarue2016}. However, here, we follow a different approach  to incorporate the common process in the major-minor LQG mean field game framework. First, in Section \ref{sec:MFevolution}, the evolution of the state mean field and the control mean field in the infinite population case are derived. Then, an $\mathcal{F}^{0,r}$-adapted and $\mathcal{F}^{i,r}$-adapted, $i \in \mathfrak{N}$, projections of the common process in the infinite population case are presented in Section \ref{sec:CommonProcess-InfPop}. Next in Sections \ref{sec:majorInf} and \ref{sec:minorInf}, the common process is perceived as a major agent in the major-minor LQG MFG framework. Subsequently, the major-minor LQG analysis described above is further extended where the major agent's state is extended with the mean field and the $\mc{F}^{0,r}$-adapted common process, while a minor agent's state is extended with the major agent's state, the mean field, and the $\mc{F}^{i,r}$-adapted common process. Finally, a convex analysis approach is taken in Section \ref{sec:SMmfgCC} to obtain the best response strategies which yield the infinite population Nash equilibrium and finite population $\epsilon$-Nash equilibrium.

\subsection{Mean Field Evolution} \label{sec:MFevolution}
The common process $y_t$ governed by \eqref{GSCMFGminorComDynNonEst} appears in the empirical average of the minor agents' states, i.e. $x^{(N)}_t$, as well as the empirical average of the minor agents' control actions, i.e. $u^{(N)}_t$. Hence, to attain the infinite population limit $\bar{y}_t$ of $y_t$, we introduce the state mean field $\bar{x}_t$ and the control mean field $\bar{u}_t$ as the infinite population limits of $x^{(N)}_t$ and $u^{(N)}_t$, respectively.
\subsubsection{Control Mean Field}
The empirical average of minor agents' control actions is defined as 
\begin{align}
u_t^{(N_k)} = \frac{1}{N_k} \sum_{j=1}^{N_k} u_{t}^{j,k},~~~~~k \in \mfK,  
\end{align} 
and we define the vector $u^{(N)}_t:= [u_t^{(N_1)}, u_t^{(N_2)}, ..., u_t^{(N_K)}]$, where the limit (in quadratic mean) of $u^{(N)}_t$, if it exists, is called the control mean field of the system and is denoted by $\bar{u}_t= [\bar{u}^1_t, ..., \bar{u}_t^K]$. 
For each minor agent $\mathcal{A}_i,\, i\in\mfN,$ of type $k\in\mfK$, we consider a state feedback control $u^{i,k}_t \in \mathcal{U}^{i,L}$ of the form 
\begin{align} \label{LstateFBcontrol}
u_t^{i,k} = L_1^k x^{i,k}_t + \sum_{l=1}^{K} \sum_{j=1}^{N_l} L_2^{k,l} x^{j,l}_t + L_3^k x^0_t + L_4^k y_t + m^k_t,
\end{align}
where $t\in[0,T]$, $L_1^k$, $L_2^{k,l}$, $L_3^{k}$ and $L_4^k$ are constant matrices of appropriate dimensions, and $m^k_t $ is a $\mc{F}^y_t$-measurable process. Coefficients in \eqref{LstateFBcontrol} depend only on agent type and not the individual agent. If we take the average of the control actions $u_t^{i,k}$ over the population of the agents of type  $k\in\mfK$, and hence calculate $u^{(N)}_t$, it can be shown that the limit $\bar{u}_t$ of $u^{(N)}_t$ as $N \rightarrow \infty$, i.e. the control mean field is given by
\begin{align} 
\bar{u}_t & = \bar{C} \bar{x}_t + \bar{D} x^0_t +\bar{E} \bar{y}_t +  \bar{r}_t, \label{GSCMFGctrlMF}
\end{align}
where $\bar{x}_t$, if it exists, denotes the state mean field introduced in Section \ref{sec:stateMF}, $\bar{y}_t$ denotes the limiting process associated with the common process $y_t$ as $N \rightarrow \infty$ (see Section \ref{sec:CommonProcess-InfPop}), and $\bar{r}_t$ is a $\mc{F}^y_t$-measurable process. Furthermore, the matrices in \eqref{GSCMFGctrlMF}, i.e.
\begin{align}
\bar{C} = \left[ \begin{array}{c}
\bar{C}_1\\
\vdots \\
\bar{C}_K \end{array} \right], \quad
\bar{D} = \left[ \begin{array}{c}
\bar{D}_1\\
\vdots\\
\bar{D}_K \end{array} \right], \quad
\bar{E} = \left[ \begin{array}{c}
\bar{E}_1\\
\vdots\\
\bar{E}_K \end{array} \right],  \quad 
\bar{r}_t = \left[ \begin{array}{c}
\bar{r}_t^1\\
\vdots\\
\bar{r}_t^K \end{array} \right],
\end{align}
are to be solved for using the mean field consistency equations \eqref{GSCMFGlatentCond1}-\eqref{GSCMFGlatentCond2} derived  in Section \ref{sec:SMmfgCC}.

\subsubsection{State Mean Field}\label{sec:stateMF}
Similarly, the empirical state average is defined as
\begin{align}
x_t^{(N_k)} = \frac{1}{N_k} \sum_{j=1}^{N_k} x_{t}^{j,k}, \qquad k\in\mfK,  
\end{align}  
where the vector $x^{(N)}_t:= [x_t^{(N_1)}, x_t^{(N_2)}, ..., x_t^{(N_K)}]$, and the limit (in quadratic mean) of $x^{(N)}_t$, if it exists, is called the state mean field of the system and is denoted by $\bar{x}_t= [\bar{x}^1_t, ..., \bar{x}_t^K]$. 

If we substitute \eqref{LstateFBcontrol} in (\ref{GSCMFGMinorIncomDyn}) for $i\in\mfN$, and take the average of the states of the minor agents' closed-loop systems of type $k\in\mfK$,
%\seb{what does ``minor agents' closed-loop systems of type $k\in\mfK$'' mean?}\dena{In systems and control terminology, when the control action is a function of the state, it is called closed-loop control, and when it is substituted in the dynamics, it is called closed-loop dynamics.} 
and hence calculate $x^{(N)}_t$, it can be shown that the limit $\bar{x}_t$ of $x^{(N)}_t$, i.e. the state mean field, satisfies 
\begin{align}
d\bar{x}_t = \left(\bar{A} \bar{x}_t + \bar{G} x^0_t + \bar{L} \bar{y}_t + \bar{m}_t\right) dt, \label{GSCMFGstateMF}   
\end{align}
where $\bar{y}_t$ denotes the infinite population limit of the common process $y_t$ (see Section \ref{sec:CommonProcess-InfPop}), $\bar{m}_t$ is a $\mc{F}^y_t$-measurable process, and the matrices
\begin{align}
\bar{A} = \left[ \begin{array}{c}
\bar{A}_1\\
\vdots \\
\bar{A}_K \end{array} \right], \quad
\bar{G} = \left[ \begin{array}{c}
\bar{G}_1\\
\vdots \\
\bar{G}_K \end{array} \right], \quad
\bar{L} = \left[ \begin{array}{c}
\bar{L}_1\\
\vdots \\
\bar{L}_K \end{array} \right], \quad
\bar{m}_t = \left[ \begin{array}{c}
\bar{m}^1_t\\
\vdots \\
\bar{m}^K_t \end{array} \right],
\end{align}
are again to be solved for using the mean field consistency conditions \eqref{GSCMFGlatentCond1}-\eqref{GSCMFGlatentCond2} derived in Section \ref{sec:SMmfgCC}.

 With a slight abuse of language, we refer to the mean value of the system's  mean field,
 %\seb{I am not sure what  ``system's Gaussian mean field'' means here} \dena{I think the intention of this statement is as follows. In general, the mean field referes to the distribution of a generic agent. In the LQG case the distribution becomes Gaussian and the mean value is referred to as the system mean field. } 
 given by the state process $\bar{x}_t = [\bar{x}^1_t, ..., ~\bar{x}^K_t]$, as  the system's mean field (the derivation of the state mean field equation above may be performed using the methods of  \cite{KizilkaleTAC2016}, \cite{Kizilkale2013} and \cite{Huang2010}).
 
\subsection{Common Process: Infinite Population} \label{sec:CommonProcess-InfPop}
Each agent completely observes the common process $y_t$ but has no observations on the latent Markov chain process $\Gamma_t$.  
In order to resolve the issue of the unobserved latent process $\Gamma_t$, the Wonham filtering method is used to estimate the distribution of $\Gamma_t$ based on the observations of each agent on $y_t$, i.e. $\mc{F}^y_t$. As a result of  \textit{Assumptions \ref{ass:majorObs}-\ref{ass:minorObs}}, the major agent $\mc{A}_0$, and each minor agent $\mc{A}_i,\, i \in \mathfrak{N}$, completely observe the unaffected common process $y^L_t$ given by \eqref{GSCMFGminorComDynNonEst} in the infinite population limit. Subsequently, $f(t,y_t^L, \Gamma_t)$ and $w_t$ in \eqref{GSCMFGminorComDynNonEst-unimpacted} are presented in their $\mc{F}^{y^L}_t$-adapted forms (see e.g. \cite{JaimungalPhil2016}, \cite{Wonham1964}).

Denote the transition probabilities for the continuous time Markov chain process $\Gamma$ by
\begin{align}
p_{ij} = \mathbf{P}(\Gamma_{t+h} = \gamma_j | \Gamma_t = \gamma_i ),\quad i,j\in\mfM
\end{align}
and the corresponding transition rates by $v_{ij} \geq 0$, and 
\begin{align}
v_i = \sum_{j=1,\, j\neq i}^{M} v_{ij}, \quad i \in \mathfrak{M}.
\end{align}
The posterior distribution of $\Gamma_t$ conditional on $\mc{F}_t^{y^L}$ is denoted by $\Pi = \{ \pi^j_t,\, j \in \mathfrak{M},\, t \in [0,T] \}$, where 
\begin{align}
\pi^j_t = \mathbb{E}[\mathds{1}_{\{ \Gamma_t = \gamma_j \}}| \mc{F}^{y^L}_t ], \quad j \in \mathfrak{M}, \quad t \in [0,T], 
\end{align}
with initial distribution $\{ \pi^j_0, j \in \mathfrak{M} \}$.
%\begin{remark}
%As a result of  \textit{Assumptions \ref{ass:majorObs}-\ref{ass:minorObs}}, the major agent $\mc{A}_0$, and each minor agent $\mc{A}_i,\, i \in \mathfrak{N}$, completely observe the unaffected common process $y^L_t$ given by \eqref{GSCMFGminorComDynNonEst} in the infinite population limit. Consequently \seb{this is not obvious... we know that under there exists a Nash equilibria with  the mean-field being $\mathcal{F^{y^L}}$--adapted, but is it a consequence of 3--5?} 
%\begin{align} 
%\pi^j_t = \mathbb{E}[\mathbbm{1}_{\{ \Gamma_t = \gamma_j \}}| \mc{F}^y_t ] \triangleq \mathbb{E}[\mathbbm{1}_{\{ \Gamma_t = \gamma_j \}}| \mc{F}^{y^L}_t ].
%\end{align}
%\end{remark}
\begin{lemma}[Wonham Filter]\cite{Wonham1964} \label{lem:WonhamFilter} 
If $\sigma > 0$, the posterior distribution $\Pi$ of $\Gamma_t$ solves the SDE
%\seb{some labeling seems off... $j$ appears to be free in the rhs, but it does not appear in the lhs, and $i$ is being summed over in the rhs.} \dena{This was corrected in the latest version of the paper.}
\begin{multline}
d\pi^j_t = \Big(- v_j \pi^j_t + \sum_{i=1,\, i\neq j}^{M} v_{ij} \pi^i_t \Big) dt - \sigma^{-2} \Big(\sum_{i=1}^{M} \pi^i_t \gamma_i\Big) \Big[f(t, y_t^L, \gamma_j) - \sum_{i=1}^{M} \pi^i_t \gamma_i\Big] \pi^j_t dt\\
+ \sigma^{-2} \Big[f(t,y^L_t, \gamma_j)-\sum_{i=1}^{M} \pi^i_t \gamma_i\Big] \pi^j_t dy_t^L,
\end{multline}  
$j \in \mathfrak{M}$.
\end{lemma}
\hfill $\square$
\begin{lemma}\cite{JaimungalPhil2016} \label{lem:PhilSeb}
Define the process $\widehat{w} = (\widehat{w}_t)_{t \in [0,T]}$ as 
\begin{align}
\widehat{w}_t = w_t  + \sigma^{-1} \int_0^t  (f_{\tau} - \widehat{f}_{\tau}) d\tau,
\end{align} 
where $\widehat{f} = (\widehat{f}_t)_{t \in [0,T]}$ is the $\mc{F}_t^y$-adapted process
\begin{align}
\widehat{f}_{t}=  \mb{E} [f(t,y_t^L, \Gamma_t)| \mathcal{F}_t^y]
% \end{align} 
% and is computed by 
% \begin{align}
% \widehat{f}_t  = \widehat{f}(t, y_t^L,\Pi ) 
= \sum_{j=1}^{M} \pi^j_t f(t,y_t^L,\gamma_j).
\end{align}
Then the process $\widehat{w}$ is an $\mc{F}^y$-adapted Wiener process. 
\end{lemma}
\hfill $\square$

Using \textit{Lemma \ref{lem:WonhamFilter}} and \textit{Lemma \ref{lem:PhilSeb}}, equation \eqref{GSCMFGminorComDynNonEst-unimpacted} can be rewritten as 
\begin{equation} \label{CommonLhat}
dy_t^L  = \widehat{f}_t\,dt + \sigma\,d\widehat{w}_t, 
\end{equation} 
and by substituting \eqref{CommonLhat} in \eqref{GSCMFGminorComDynNonEst}, the $\mc{F}^{0,r}$-adapted dynamics of the common process for the infinite population case, i.e. $\bar{y}_t$, is given by 
\begin{align} \label{GSCMFGComDyn}
d\bar{y}_t = [\widehat{f}_{t}  + F^{\pi} \bar{u}_t + F_0 u^0_t + H^{\pi} \bar{x}_t + H_0 x^0_t]\,dt + \sigma \,d\widehat{w}_t, 
\end{align}
where the average terms $x^{(N)}_t$ and $u^{(N)}_t$ in \eqref{GSCMFGminorComDynNonEst} have been replaced with their limits as $N \rightarrow \infty$, i.e.  the state mean field $\bar{x}_t$ and the control mean field $\bar{u}_t$, respectively. Moreover, $F^{\pi} = \pi \otimes F$ and $H^{\pi} = \pi \otimes H$, where $\otimes$ denotes the Kronecker product of the corresponding matrices. 
\begin{remark}
As the state and the control action of each individual minor agent $\mc{A}_i,\, i \in \mathfrak{N}$, do not affect the infinite population evolution of the common process, i.e.  $\bar{y}_t$, the $\mc{F}^{i,r}$-adapted and $\mc{F}^{0,r}$-adapted projection of the common  process $\bar{y}_t$ in the infinite population limit  are identical and given by \eqref{GSCMFGComDyn}.
\end{remark}

\subsection{Major Agent's Regulation Problem : Infinite Population} \label{sec:majorInf}
To begin, we extend the major agent's state $x^0_t$  with the state mean field $\bar{x}_t$ and the infinite population common process $\bar{y}_t$ to form the major agent's extended state $X_t^{0} = [ (\bar{y}_t)^\intercal, \; (x^0_t)^\intercal, \; (\bar{x}_t)^\intercal ]^\intercal$ and satisfies the SDE 
\begin{align} dX^{0}_t = \mb{A}_0 X^{0}_t dt + \mb{B}_0 u^0_t dt + \mb{M}_t^0 dt + \Sigma_0 dW^0_t.
\label{GSCMFGlatentMajorExtDyn} 
\end{align}
By substituting \eqref{GSCMFGctrlMF} into \eqref{GSCMFGComDyn}, the matrices in the extended major agent's dynamics \eqref{GSCMFGlatentMajorExtDyn} are given by 
\begin{gather}
\mb{A}_0 = \left[ \begin{array}{ccc}
F^{\pi} \bar{E} & F^{\pi} \bar{D} + H_0 &F^{\pi} \bar{C} + H^{\pi}  \\
0_{n \times n} & A_0 & 0_{n \times nK}\\
\bar{L} & \bar{G} & \bar{A}
\end{array} \right], \quad
\mb{B}_0 = \left[ \begin{array}{c}
F_0\\
B_0 \\
0_{nK \times m}
\end{array} \right],
\quad
\mb{M}^0_t = \left[ \begin{array}{c}
\widehat{f}_t+ F^{\pi} \bar{r}_t\\
b_0(t)\\
\bar{m}_t
\end{array} \right], 
\\
\Sigma_0 = \left[ \begin{array}{ccc}
 \sigma & 0_{n \times r} & 0_{n \times rK} \\
0_{n \times r} & \sigma_0 & 0_{n \times rK} \\
0_{nK \times r} & 0_{nK \times r} & 0_{nK \times rK} 
 \end{array} \right], \quad 
W^0_t = \left[ \begin{array}{c}
\widehat{w}_t\\
w^0_t\\
0_{rK \times 1}
\end{array} \right].\nonumber
\end{gather}
Next, the major agent's extended cost functional is 
\begin{equation} \label{GSCMFGlatent-majorcostExt}
J_0^{ex}(u^0) = \frac{1}{2}\mb{E} \bigg [  X_T^{0\intercal} \mb{G}_0 X_T^{0} + \int_0^T \Big \{X_s^{0\intercal}\mb{Q}_0 X_s^{0} + 2\,X_s^{0\intercal} \mb{N}_0 u_s^0 + u^{0\intercal}_s R_0 u^0_s \Big \} ds \bigg], 
\end{equation}
where the corresponding weight matrices are given by 
\begin{gather}
\mathbb{G}_0 = \left[I_{2n}, 0_{2n \times nK}\right ]^\intercal G_0 \left [I_{2n}, 0_{2n \times nK}\right],\\ \mathbb{Q}_0 = \left[I_{2n}, 0_{2n \times nK}\right ]^\intercal Q_0 \left[I_{2n}, 0_{2n \times nK}\right], \\
\mathbb{N}_0 = \left[ \begin{array}{c} 
N_{0}\\
0_{nK \times m}  \end{array} \right].
\end{gather}
The minimization of the extended cost functional \eqref{GSCMFGlatent-majorcostExt} subject to the extended dynamics \eqref{GSCMFGlatentMajorExtDyn} constitutes a stochastic optimal control problem for the major agent in the infinite population limit. Our main result  \textit{Theorem \ref{thm}} is proven in Section \ref{sec:SMmfgCC}, and it implies that the major agent's optimal control action is given by
 \begin{align}\label{GSCMFGLatentSLQconvControl-major}
 u^{0,*}_t = - R_0^{-1} \Big [ \mathbb{N}_0^\intercal X^{0}_t +  \mathbb{B}_0^\intercal \big(\Pi_0(t) X_t^{0} + s^0_t \big) \Big],
 \end{align}
 where $\Pi_0(t)$ is the solution %\seb{do we have uniqueness?} 
 to the ODE
 \begin{subequations}
 \label{GSCMFGlatentLQGriccati-major}
\begin{align}
-\dot{\Pi}_0 = \Pi_0 \mb{A}_0 + \mb{A}_0^\intercal \Pi_0 - (\mb{B}_0^\intercal \Pi_0 + \mb{N}_0^\intercal)^\intercal R_0^{-1}( \mb{B}_0^\intercal \Pi_0 + \mb{N}_0^\intercal ) + \mb{Q}_0  &, \\
\Pi_0(T) &= \mb{G}_0,
\end{align} 
\end{subequations}
and $s^0_t$ is the solution %\seb{do we have uniqueness?}
to the  FBSDE %\seb{what is $q_t^0$, have all other objects been defined or referenced? Can we emphasize where is the forward part... I presume some of the coefficients such as $\hat{M}_t$... but it could be useful to point that out.}\dena{I defined $q^0_t$ and made a forward reference to the exact definition, all other coefficients are defined beforehand, the forward part is $\mb{M}^0_t$ which I mentioned it just under the equation. }
\begin{subequations}
\begin{align}
%\begin{split}
 d{s}^0_t = - \left( [(\mb{A}_0- \mb{B}_0R_0^{-1} \mb{N}_0)^\intercal- \Pi_0 \mb{B}_0 R_0^{-1} \mb{B}_0^\intercal ] s^0_t +  \Pi_0 \mb{M}^0_t\right)dt %& \\
 + (q^0_t - \Pi_0 \Sigma_0)dW^0_t,% &,
%\end{split}
 \\
s^0_T &= 0,
\end{align} 
\label{GSCMFGLatentLQGoffset-major}
\end{subequations}
where $\mb{M}^0_t$ is generated through forward propagation and $q^0_t$ is an $\mathcal{F}^{0,r}$-adapted process defined under equation \eqref{GSCMFGLatentSLQconvdiff2-major}. The Riccati equation \eqref{GSCMFGlatentLQGriccati-major} and the offset equation \eqref{GSCMFGLatentLQGoffset-major} shall be derived in Section \ref{sec:SMmfgCC}.

We use the notation
\begin{gather}
g^0_t = g^0(t, y^L_t, s^0_t) =  [(\mb{A}_0- \mb{B}_0R_0^{-1} \mb{N}_0)^\intercal- \Pi_0 \mb{B}_0 R_0^{-1} \mb{B}_0^\intercal ] s^0_t +  \Pi_0 \mb{M}^0_t,\\
\xi^0_t = [y^L_t, s^0_t]^\intercal, \\
 A^0(t, \xi^0) = \begin{bmatrix}
 - H^\intercal_0 g^0\\
 H_0 \hat{f}
 \end{bmatrix}(t, y^L, s^0), 
\end{gather}
for a fixed $(2n + nK) \times n$ full-rank matrix $H_0$.%\dena{I added theses assumptions for the existence and uniqueness of the solution to the BFSDE.}
%We use the notation $\xi^0$ for the vector $\xi^0 = [y^L, s^0, q^0]^\intercal$, and for a fixed $n \times (2n + nK)$ full-rank matrix $H_0$, we denote 
%\begin{equation}
% A^0(t, \xi^0) = \begin{bmatrix}
% - H^\intercal_0 g^0\\
% H_0 \hat{f}
% \end{bmatrix}(t, y^L, s^0, q^0).
%\end{equation}
\begin{assumption}[Measurability Condition] The processes $\hat{f}_t$ and $g^0_t$ are, respectively, $\mdsR^n$-valued and $\mdsR^{2n+nK}$-valued progressively measurable processes on $[0,T]$ such that (see \cite{CarmonaBSDEsBook2016})
\begin{gather}
\mb{E}\left[\int_0^T \Vert \hat{f}_t \Vert^2 \, dt\right] < \infty,\\
\mb{E}\left[\int_0^T \Vert g_t^0 \Vert^2 \, dt\right] < \infty.
\end{gather}

\end{assumption}
\begin{assumption}[Growth Condition] There exists a constant $c>0$ such that for each $(t, \omega, y^L, s^0) \in [0,T] \times \Omega \times \mdsR^d \times \mdsR^{2n+nK}$, we have (see \cite{CarmonaBSDEsBook2016})
\begin{gather}
  \Vert\hat{f}(t, y^L) \Vert \leq c(1+ \Vert y^L \Vert),\\
  \Vert g^0(t, y^L, s^0) \Vert \leq c(1 + \Vert y^L \Vert + \Vert s^0 \Vert ).
\end{gather}

\end{assumption}

\begin{assumption}[Lipschitz Condition] There exists a constant $c > 0$ such that $\forall t \in [0, T]$, $\forall \omega \in \Omega$, $\forall y^L, {y^L}^\prime \in \mdsR^d$, and $\forall s^0, {s^0}^{\prime} \in \mdsR^{2n+nK}$, we have (see \cite{CarmonaBSDEsBook2016})
\begin{gather}
  \Vert \hat{f}(t, y^L) - \hat{f}(t, {y^L}^\prime)\Vert \leq c(1+ \Vert y^L - {y^L}^\prime \Vert),\\
  \Vert g^0(t, y^L, s^0)- g^0(t, {y^L}^\prime, {s^0}^\prime) \Vert \leq c(1 + \Vert y^L-{y^L}^\prime \Vert + \Vert s^0-{s^0}^\prime \Vert).
\end{gather}
\end{assumption}

\begin{assumption}[Monotonicity Condition] There exist $\beta_1$ and $\beta_2$ such that  $\forall t \in [0, T]$, $\forall \omega \in \Omega$, $\forall y^L, {\tilde{y}^L} \in \mdsR^d$, and  $\forall s^0, {\tilde{s}^0} \in \mdsR^{2n+nK}$, we have (see \cite{CarmonaBSDEsBook2016})
\begin{equation}
    [A^0(t, \xi^0) - A^0(t, \tilde{\xi}^0)].(\xi^0 - \tilde{\xi}^0) \leq - \beta_1 \Vert H_0(y^L - \tilde{y}^L) \Vert ^2 -\beta_2 \Vert H^T_0(s^0-\tilde{s}^0) \Vert^2,
\end{equation}
where $\beta_1 > 0$, $\beta_2 \geq 0$. 
\end{assumption}
Finally, the closed-loop dynamics of the major agent $\mc{A}_0$ when the control action \eqref{GSCMFGLatentSLQconvControl-major} is substituted in \eqref{GSCMFGMajorIncomDyn} satisfies the SDE
\begin{equation}
\label{majorDycsClosed}
dX_t^0 = \left ( A_0X_t^0 - B_0 R_0^{-1} [ \mathbb{N}_0^\intercal X^{0}_t +  \mathbb{B}_0^\intercal (\Pi_0(t) X_t^{0} + s^0_t ) ] + b_0(t) \right)dt + \sigma_0 \,dw^0_t.  
\end{equation}

\subsection{Minor Agent's Regulation Problem: Infinite Population}  \label{sec:minorInf}
For  minor agent $\mc{A}_i$'s, $i \in \mathfrak{N}$, we extend the state with the infinite population common process $\bar{y}_t$, the major agent's state $x^0_t$, and the state mean field $\bar{x}_t$ to form the minor agent's extended state $X^{i}_t = [(x^i_t)^\intercal,\, (\bar{y}_t)^\intercal,\, (x^0_t)^\intercal, (\bar{x}_t)^\intercal]^\intercal$ which satisfies  
\begin{equation}\label{GSCMFGlatant-minorExtDyn}
dX^{i}_t = \left(\mb{A}_k X^{i}_t  + \mb{B}_k u_i  + \mb{M}_t^k\right) dt + \Sigma_k \,dW^i_t.
\end{equation}
To attain the extended matrices in \eqref{GSCMFGlatant-minorExtDyn}, the joint dynamics of (i) minor agent $\mc{A}_i$'s system given by \eqref{GSCMFGMinorIncomDyn}, (ii) the common process $\bar{y}_t$ given by \eqref{GSCMFGComDyn} where \eqref{GSCMFGctrlMF} and \eqref{GSCMFGLatentSLQconvControl-major} are substituted , (iii) the major agent $\mc{A}_0$'s closed loop system  given by \eqref{majorDycsClosed}, and (iv) the state mean field $\bar{x}_t$ given by \eqref{GSCMFGstateMF} are utilized which results in 
\begin{gather}
\mb{A}_k = \left[ \begin{array}{cc}
A_k & 0_{n \times (2n+nK)} \\
0_{(2n+nK) \times n} & \mb{A}_0 - \mb{B}_0 R_0^{-1} \mb{N}_0^\intercal -\mb{B}_0 R_0^{-1} \mb{B}_0^\intercal \Pi_0 
\end{array} \right], \quad
\mb{B}_k = \left[ \begin{array}{c}
B_k \\
0_{(2n+nK) \times m}
\end{array} \right],\\
\mb{M}^k_t = \left[ \begin{array}{c}
b_k(t)\\
\mb{M}^0_t - \mb{B}_0 R_0^{-1} \mb{B}_0^\intercal s_0 (t)
\end{array} \right],\quad
\Sigma_k = \left[ \begin{array}{cc}
\sigma_k & 0_{n \times (2r+rK)} \\
 0_{(2n+nK) \times r} & \Sigma_0  
\end{array} \right], \quad 
W^i_t = \left[ \begin{array}{c}
w^i_t,\\
W^0_t
\end{array} \right].\nonumber
\end{gather}
 
Next, the minor agent $\mc{A}_i$'s extended cost functional is formed as 
\begin{equation}\label{GSCMFGlatent-minorcostExt}
J_i^{ex}(u^i) = \frac{1}{2}\mb{E} \bigg [  X_T^{i\intercal} \mb{G}_k X_T^{i} + \int_0^T \Big \{X_s^{i\intercal} \mb{Q}_k X_s^{i} + 2\,X_s^{i\intercal} \mb{N}_k u_s^i + u^{i\intercal}_s R_k u^i_s \Big \} ds \bigg], 
\end{equation}
where the corresponding weight matrices are given by  
\begin{gather}
\mathbb{G}_k = \left[I_{2n}, 0_{2n \times (n+nK)}\right ]^\intercal G_k \left [I_{2n}, 0_{2n \times (n+nK)}\right] ,\nonumber \\ \mathbb{Q}_k = \left[I_{2n}, 0_{2n \times (n+nK)}\right ]^\intercal Q_k \left[I_{2n}, 0_{2n \times (n+nK)}\right], \nonumber\\
\mathbb{N}_k = \left[ \begin{array}{c} 
N_{k}\\
0_{(n+nK) \times m}  \end{array} \right].
\end{gather}
The dynamics \eqref{GSCMFGlatant-minorExtDyn} together with the cost functional \eqref{GSCMFGlatent-minorcostExt} constitute a stochastic optimal control problem for  minor agent $\mc{A}_i,\, i \in \mathfrak{N}$, in the infinite population limit. Once again, according to \textit{Theorem \ref{thm}}, the minor agent $\mc{A}_i$'s optimal control action for the infinite population case is given by  
 \begin{equation}\label{GSCMFGLatentSLQconvControl}
 u^{i,*}_t = - R_k^{-1} \Big [ \mathbb{N}_k^\intercal X^{i}_t +  \mathbb{B}_k^\intercal \big(\Pi_k(t) X_t^{i} + s^{i,k}_t \big) \Big],
 \end{equation} 
 where $\{\Pi_k(t), \, k \in \mfK\}$, is the solutions to the following coupled %\seb{do we have uniqueness?}
 deterministic Riccati equation 
 \begin{subequations}
 \label{GSCMFGlatentLQGriccati}
\begin{align}
-\dot{\Pi}_k = \Pi_k \mb{A}_k + \mb{A}_k^\intercal \Pi_k - (\mb{B}_k^\intercal \Pi_k + \mb{N}_k^\intercal)^\intercal R_k^{-1}( \mb{B}_k^\intercal \Pi_k + \mb{N}_k^\intercal ) + \mb{Q}_k  &, \\
\Pi_k(T) &= \mb{G}_k,
\end{align} 
\end{subequations}
and $s^{i,k}_t,\, k \in \mfK$, is the solution %\seb{do we have uniqueness?}
to the following FBSDE %\seb{what is $q_t^i$, have all other objects been defined or referenced? Can we emphasize where is the forward part... I presume some of the coefficients such as $\hat{M}_t$... but it could be useful to point that out.} \dena{I defined $q^i_t$ and mentioned the forward part, all other coefficients are defined previously.}
\begin{subequations}
\label{GSCMFGLatentLQGoffset}
\begin{align} 
%\begin{split}
 ds^{i,k}_t =- \Big( \big[(\mb{A}_k- \mb{B}_kR_k^{-1} \mb{N}_k)^\intercal- \Pi_k \mb{B}_k R_k^{-1} \mb{B}_k^\intercal \big] s^{i,k}_t
 + \Pi_k \mb{M}^k_t \Big )dt %&
 %\\
 + (q^i_t - \Pi_k \Sigma_k)d\,W^i_t,% &,
%\end{split}
\\
s^{i,k}_T &= 0,
\end{align}
\end{subequations}
where $\mb{M}^k_t$ is generated through forward propagation, and $q^i_t$ is an $\mc{F}^{i,r}$-adapted process defined under equation \eqref{GSCMFGLatentSLQconvdiff2}. The complete derivation of \eqref{GSCMFGlatentLQGriccati}-\eqref{GSCMFGLatentLQGoffset} will be discussed in Section \ref{sec:SMmfgCC}. 

We use the notation \begin{gather}
g^{i,k}_t = g^{i,k}(y^L_t, s^i_t) =  [(\mb{A}_k- \mb{B}_kR_k^{-1} \mb{N}_k)^\intercal- \Pi_k \mb{B}_k R_k^{-1} \mb{B}_k^\intercal ] s^{i}_t +  \Pi_k \mb{M}^k_t,\\
\xi^i_t = [y^L_t, s^i_t]^\intercal, \\
 A^{i,k}(t, \xi^i) = \begin{bmatrix}
 - H^\intercal_k g^{i,k}\\
 H_k \hat{f}
 \end{bmatrix}(t, y^L, s^i), 
\end{gather}
for a fixed $(3n + nK) \times n$ full-rank matrix $H_k$. %\dena{I added the following assumptions for the existence and uniqueness of the solution to the BFSDE.}

\begin{assumption}[Measurability Condition] The processes 
 $\hat{f}_t$ and $g^0_t$ are, respectively, $\mdsR^n$-valued and $\mdsR^{3n+nK}$-valued progressively measurable processes on $[0,T]$ such that (see \cite{CarmonaBSDEsBook2016})
\begin{gather}
\mb{E}\left[\int_0^T \Vert \hat{f}_t \Vert^2 \, dt\right] < \infty,\\
\mb{E}\left[\int_0^T \Vert g_t^{i,k} \Vert^2 \, dt\right] < \infty.
\end{gather}

\end{assumption}
\begin{assumption}[Growth Condition] There exists a constant $c>0$ such that for each $(t, \omega, y^L, s^i)\in [0,T] \times \Omega \times \mdsR^d \times \mdsR^{3n+nK}$, $\forall i \in \mfN $, and $k \in \mfK$, we have (see \cite{CarmonaBSDEsBook2016})
\begin{gather}
 \Vert \hat{f}(y^L) \Vert \leq c(1+ \Vert y^L \Vert),\\
  \Vert g^{i,k}(y^L, s^i) \Vert \leq c(1 + \Vert y^L \Vert + \Vert s^i \Vert).
\end{gather}

\end{assumption}

\begin{assumption}[Lipschitz Condition] There exists a constant $c > 0$ such that $\forall t \in [0, T]$, $\forall \omega \in \Omega$, $\forall y^L, {y^L}^\prime \in \mdsR^d$, $\forall s^i, {s^i}^{\prime} \in \mdsR^{3n+nK}$, $\forall i \in \mfN $, and $k \in \mfK$, we have (see \cite{CarmonaBSDEsBook2016})
\begin{gather}
  \Vert \hat{f}(y^L) - \hat{f}({y^L}^\prime)\Vert \leq c(1+ \Vert y^L - {y^L}^\prime \Vert),\\
  \Vert g^{i,k}(y^L, s^i)- g^{i,k}({y^L}^\prime, {s^i}^\prime) \Vert \leq c(1 + \Vert y^L-{y^L}^\prime \Vert + \Vert s^i-{s^i}^\prime \Vert ).
\end{gather}
\end{assumption}
 
 \begin{assumption}[Monotonicity Condition] There exist $\beta_1$ and $\beta_2$ such that $\forall t \in [0, T]$, $\forall \omega \in \Omega$, $\forall y^L, {\tilde{y}^L} \in \mdsR^d$, and $\forall s^i, {\tilde{s}^i} \in \mdsR^{3n+nK}$, $\forall i \in \mfN $, $k \in \mfK$
We have  (see \cite{CarmonaBSDEsBook2016})
\begin{equation}
    [A^{i,k}(t, \xi^i) - A^{i,k}(t, \tilde{\xi}^i)].(\xi^i - \tilde{\xi}^i) \leq - \beta_1 \Vert H_k(y^L - \tilde{y}^L)\Vert^2 -\beta_2 \Vert H^T_k(s^i-\tilde{s}^i) \Vert^2,
\end{equation}
where $\beta_1 > 0$, $\beta_2 \geq 0$. 
\end{assumption}
Finally, control action \eqref{GSCMFGLatentSLQconvControl} %\seb{incorrect ref... or multiply defined ref}\dena{corrected.} 
is substituted in \eqref{GSCMFGMinorIncomDyn} which gives minor agent $\mc{A}_i$'s, $i \in \mathfrak{N}$, closed loop system as 
\begin{equation}\label{minorClosedloop}
dX^i_t = \Big ( A_k X^i_t - B_k R_k^{-1} \Big [ \mb{N}_k^\intercal X^{i}_t + \mb{B}_k^\intercal \big (\Pi_k X_t^{i} + s^{i,k}_t\big) \Big ] + b_k \Big ) dt+ \sigma_k dw^i_t.
\end{equation}
\begin{remark}
When there are no latent processes, i.e. $y^L_t =0,\, t \in [0,T]$, the diffusion terms %\seb{can we say which  terms precisely are zero?} \dena{I added equations (45) and (46). In this case $q^i_t$ is a deterministic function of time which can be specified wrt the type of the agent.} 
of \eqref{GSCMFGLatentLQGoffset-major} and \eqref{GSCMFGLatentLQGoffset} become zero as we have the deterministic functions of time $q^0$ and $q^k$ given by
\begin{gather}
  q^0 = \Pi_0 \Sigma_0,\\
  q^k = \Pi_k \Sigma_k,
\end{gather}
which for minor agents are specified with respect to their type. Hence equations \eqref{GSCMFGLatentLQGoffset-major} and \eqref{GSCMFGLatentLQGoffset} reduce to the deterministic offset equations of classical major-minor LQG mean field games in \cite{Huang2010}. 

\end{remark}
\subsection{ Nash and $\epsilon$-Nash Equilibria } \label{sec:SMmfgCC}
To derive the mean field consistency equations which specify the matrices in the control and state mean field equations \eqref{GSCMFGctrlMF} and \eqref{GSCMFGstateMF}, respectively, the closed loop system \eqref{minorClosedloop} of minor agent $\mc{A}_i$ may be written out explicitly as
\begin{multline} \label{GSCMFGLatentMinorsCtrld}
dx^i_t = \Big ( A_k x^i_t - B_k R_k^{-1} \big (\mb{N}_k^\intercal + \mb{B}_k^\intercal \Pi_k \big) \big[ x^{i\intercal}_t, \bar{y}_t^\intercal, x^{0\intercal}_t, \bar{x}_t^\intercal\big]^\intercal \\
- B_k R_k^{-1} \mb{B}_k^\intercal s^{i,k}_t + b_k \Big ) dt + \sigma_k \,dw^i_t,
\end{multline}
where $i \in \mathfrak{N}$, $ k \in \mfK$.

Next, define the block matrices 
\begin{gather}
\Pi_k = \left[ \begin{array}{cccc}
\Pi_{k,11} & \Pi_{k,12} & \Pi_{k,13} & \Pi_{k,14} \\
\Pi_{k,21} & \Pi_{k,22} & \Pi_{k,23} & \Pi_{k,24}\\
\Pi_{k,31} & \Pi_{k,32} & \Pi_{k,33} & \Pi_{k,34}\\
\Pi_{k,41} & \Pi_{k,42} & \Pi_{k,43} & \Pi_{k,44}
\end{array} \right], \quad 
\mb{N}_k = \left[ \begin{array}{c}
\mb{N}_{k,1}\\
\mb{N}_{k,2}\\
\mb{N}_{k,3}\\
\mb{N}_{k,4}
\end{array} \right],
\\ 
\mathbf{e}_k = \left[ 0_{n\times n}, ... ,  0_{n\times n}, I_n,  0_{n\times n}, ...,  0_{n\times n}\right ],
\end{gather}
where $\Pi_{k,11}$, $\Pi_{k,22}$, $\Pi_{k,33} \in \mdsR^{n\times n}$, $\Pi_{k, 44} \in \mdsR^{nK \times nK} $; $\mb{N}_{k,1}, \mb{N}_{k,2}, \mb{N}_{k,3} \in \mdsR^{n \times m}$, $\mb{N}_{k,4} \in \mdsR^{nK \times m}$; and $\mathbf{e}_k \in \mdsR^{n \times nK}$, denotes a matrix which has the identity matrix $I_n$ in its $k$th block and zero matrix $0_{n \times n}$ in other $(K-1)$ remaining blocks, for all $k \in\mfK$.

Now, taking the average of \eqref{GSCMFGLatentMinorsCtrld} over $N_k$ minor agents of type $k \in \mfK$, and then its limit as the number $N_k$ of agents within the subpopulation $k$ goes to infinity (i.e. $N_k \rightarrow \infty$), we find
\begin{multline} \label{GSCMFGLatent-limitAveMinors}
d\bar{x}_t^k = \Bigg\{\Big[ \big ( A_k - B_k R_k^{-1} (\mb{N}_{k,1}^\intercal + B_k^\intercal \Pi_{k, 11}) \big) \mathbf{e}_k
 - B_k R_k^{-1}(\mb{N}_{k,4}^\intercal+ B_k^\intercal \Pi_{k, 14}) \Big]  \bar{x}_t  \\
-B_k R_k^{-1} (\mb{N}_{k,3}^\intercal+B_k^\intercal \Pi_{k, 13} )x^0_t -B_k R_k^{-1} (\mb{N}_{k,2}^\intercal+B_k^\intercal \Pi_{k,12}) \bar{y}_t  +(b_k -B_k R_k^{-1} \mb{B}_k^\intercal \bar{s}^k_t)\Bigg\}dt, \quad \text{q.m.}
\end{multline}
In \eqref{GSCMFGLatent-limitAveMinors}, $\bar{s}^k_t$ is obtained by taking the average and then the limit of \eqref{GSCMFGLatentLQGoffset} over the subpopulation $k \in \mfK$ as $N_k \rightarrow \infty$, and solves
\begin{subequations}
\begin{align}
%\begin{split}
  d{\bar{s}}^k_t 
  =- \Big ( \big[(\mb{A}_k- \mb{B}_kR_k^{-1} \mb{N}_k^\intercal)^\intercal- \Pi_k \mb{B}_k R_k^{-1} \mb{B}_k^\intercal \big] \bar{s}^k_t + \Pi_k \mb{M}^k_t \Big ) dt
 % &
 % \\ 
  + (\bar{q}_t -\Pi_k \Sigma_k)d\bar{W}_t,
%\end{split}
\\
\bar{s}^k_T &= 0,
\end{align}
\end{subequations}
where $\bar{W}_t = \left[ \begin{array}{c}
0_{r\times 1},\\
W_t^{0}
\end{array} \right]$, since $\lim_{N_k \rightarrow \infty}\frac{1}{N_k}\sum_{i=1}^{N_k} w^i_t=0$; and  $\bar{q}_t=e^{-\mb{A}_k^\intercal t} \bar{Z}_t$, where $\bar{Z}_t$ is an $\mc{F}^{0,r}_t$-adapted process 
%\seb{how do we see this last claim? where is $\bar{q}_t$ defined?}\dena{I added a couple of lines here to make it clear, but it needs forward referencing again.} 
stemming from  the martingale representation of $\bar{M}_t$, which itself is the limiting average
 \begin{equation}
 \frac{1}{N} \sum_{i=1}^N M^i_t \stackrel{\text{q.m.}}{\longrightarrow} \bar{M}_t~,
\end{equation}
where $M^i_t$ is given by \eqref{martingaleMinor}.
% satisfying 
% \begin{equation}\label{GSCMFGLatentMrtnglRepLimit}
%  \bar{M}_t = \bar{M}_0 + \int_0^t \bar{Z}_s d\bar{W}_s,
%  \end{equation}
%  for the martingale $\bar{M}_t$ defined as 
%  \begin{equation}
%  \frac{1}{N} \sum_{i=1}^N M^i_t \stackrel{\text{q.m.}}{\longrightarrow} \bar{M}_t,~ 
% \end{equation}
% where $M^i_t$ is given by \eqref{martingaleMinor}.

Then, equating \eqref{GSCMFGLatent-limitAveMinors} with \eqref{GSCMFGstateMF} results in the following sets of equations.%\seb{do we really need braces around these equations? needs better formatting to fit on the lines.}\dena{I put braces because the first set only depends on the solutions of the Riccati equations and the second set only depends on the solutions of the offset equations, but I am open to use better presentations. I will reformat equations in the main template.}
\begin{align}
\begin{cases}
& -\dot{\Pi}_0 = \Pi_0 \mb{A}_0 + \mb{A}_0^\intercal \Pi_0 - (\mb{N}_0^\intercal+\mb{B}_0^\intercal \Pi_0 )^\intercal R_0^{-1}(\mb{N}_0^\intercal + \mb{B}_0^\intercal \Pi_0) + \mb{Q}_0, \quad 
\Pi_0(T) = \mb{G}_0,  \allowdisplaybreaks\\
&- \dot{\Pi}_k = \Pi_k \mb{A}_k + \mb{A}_k^\intercal \Pi_k - (\mb{N}_k^\intercal+\mb{B}_k^\intercal \Pi_k )^\intercal R_k^{-1}(\mb{N}_k^\intercal + \mb{B}_k^\intercal \Pi_k) + \mb{Q}_k  ,\quad \Pi_k(T) = \mb{G}_k,\allowdisplaybreaks \\
&\bar{C}_k= -R_k^{-1} (\mb{N}_{k,1}^\intercal + B_k^\intercal \Pi_{k, 11}) \mathbf{e}_k - R_k^{-1} (\mb{N}_{k,4}^\intercal + B_k^\intercal \Pi_{k, 14}),\allowdisplaybreaks \\
&\bar{D}_k = -R_k^{-1} (\mb{N}_{k,3}^\intercal +B_k^\intercal \Pi_{k, 13}),\allowdisplaybreaks \\
&\bar{E}_k = -R_k^{-1} (\mb{N}_{k,2}^\intercal + B_k^\intercal \Pi_{k,12} ),\allowdisplaybreaks\\
&\bar{A}_k = A_k \mathbf{e}_k + B_k \bar{C}_k,\allowdisplaybreaks \\
%[ A_k - B_k R_k^{-1} (\mb{N}_{k,1} + B_k^T \Pi_{k, 11}) ] \textbf{e}_k - B_k R_k^{-1} B_k^T \Pi_{k, 14},\\
&\bar{G}_k = B_k\bar{D}_k, \allowdisplaybreaks \\%-B_k R_k^{-1} B_k^T \Pi_{k, 13},\\
&\bar{L}_k = B_k \bar{E}_k, \label{GSCMFGlatentCond1}
%-B_k R_k^{-1} (B_k^T \Pi_{k,12} + \mb{N}_{k,2}), \label{GSCMFGlatentCond1}
\end{cases}
\end{align} 
\begin{align} \begin{cases}
 d{s}^0_t = -\Big ( \big [(\mb{A}_0- \mb{B}_0R_0^{-1} \mb{N}_0^\intercal)^\intercal- \Pi_0 \mb{B}_0 R_0^{-1} \mb{B}_0^\intercal \big] s^0_t+  \Pi_0 \mb{M}^0_t \Big) dt  + (q^0_t - \Pi_0 \Sigma_0)dW^0_t,\quad s^0_T = 0, \\
  d{\bar{s}}^k_t = -\Big ( \big[(\mb{A}_k- \mb{B}_kR_k^{-1} \mb{N}_k^\intercal)^\intercal- \Pi_k \mb{B}_k R_k^{-1} \mb{B}_k^\intercal \big]\bar{s}^k_t + \Pi_k \mb{M}^k_t \Big ) dt + (\bar{q}_t - \Pi_k \Sigma_k )d\bar{W}_t, \quad
\bar{s}^k_T = 0,\\
\bar{r}^k_t =-R_k^{-1} \mb{B}_k^\intercal \bar{s}^k_t,  \\
\bar{m}^k_t = B_k \bar{r}^k_t + b_k. \label{GSCMFGlatentCond2}
 %-B_k R_k^{-1} \mb{B}_k^T s_k + b_k.\label{GSCMFGlatentCond2}
\end{cases}
\end{align} 
Equations \eqref{GSCMFGlatentCond1}-\eqref{GSCMFGlatentCond2} are called the mean field consistency equations (see \cite{Huang2010}) from which the matrices in \eqref{GSCMFGctrlMF} and \eqref{GSCMFGstateMF} can be calculated.

Let us define the following matrix 
\begin{equation}
% M_1 = \left[ \begin{array}{ccc}
% A_1-B_1R_1^{-1}(\mb{N}_{1,1}^T+B_1^T\Pi_{1,11}) & & 0 \\
%  & \ddots & \\
% 0 & & A_K-B_KR_K^{-1}(\mb{N}_{K,1}^T+B_K^T\Pi_{K,11})
% \end{array}\right],
M_1= \text{diag}(A_k-B_kR_k^{-1}(\mb{N}_{k,1}^\intercal+B_k^\intercal \Pi_{k,11}))_{k\in\mfK}.
\end{equation}
%\begin{gather}
%M_1^\prime = \text{diag}(-\pi_kFR_k^{-1}(\mb{N}_{k,1}^T+B_k^T\Pi_{k,11}))_{k\in\mfK},
% M_1^\prime = \left[ \begin{array}{ccc}
% -\pi_1FR_1^{-1}(\mb{N}_{1,1}^T+B_1^T\Pi_{1,11}) & & 0 \\
%  & \ddots & \\
%  0 & & -\pi_KFR_K^{-1}(\mb{N}_{K,1}^T+B_K^T\Pi_{K,11})
% \end{array}\right],
%\nonumber 
%\allowdisplaybreaks
%\\ 
%M_2 = \left[ \begin{array}{c}
%-B_1R_1^{-1}(\mb{N}^T_{1,4}+B_1^T\Pi_{1,14})\\
%\vdots\\
%-B_KR_K^{-1}(\mb{N}^T_{K,4}+B_K^T\Pi_{K,14})
%\end{array}\right], \nonumber
%\\ 
%M_2^{\prime} = \left[ \begin{array}{c}
%-\pi_1FR_1^{-1}(\mb{N}^T_{1,4}+B_1^T\Pi_{1,14})\\
%\vdots\\
%-\pi_KFR_K^{-1}(\mb{N}^T_{K,4}+B_K^T\Pi_{K,14})
%\end{array}\right],\nonumber \allowdisplaybreaks\\
%M_3 = \left[\begin{array}{ccccc}
%F^{\pi}\bar{E} & F^{\pi}\bar{D}+H_0 & F^{\pi} \bar{C} + H^{\pi} & 0_{n \times nK} & 0_{n \times n}\\
%0_{n \times n} & A_0 & 0_{n \times nK} & 0_{n \times nK} & 0_{n \times n}\\
%\bar{L} & \bar{G} & \bar{A} & 0_{nK \times nK} & 0_{nK \times n}\\
%\bar{L} & \bar{G} & M_2 & M_1 & 0_{nK \times n}\\
%F^{\pi}\bar{E} & F^{\pi}\bar{D}+H_0 & M_2^{\prime} & M_1^{\prime} & 0_{n \times n} 
%\end{array}\right],\nonumber\allowdisplaybreaks\\
%L_{0,H} = Q_0^{\frac{1}{2}}\left[\begin{array}{ccccc}
%0_{n\times n} & I_{n} & 0_{n \times nK} & 0_{n \times nK} & 0_{n \times n}\\
%0_{n \times n} & 0_{n \times n} & 0_{n+nK} & 0_{n+nK} & I_n
%\end{array}\right],\nonumber \allowdisplaybreaks\\
%L_a = Q_0^{\frac{1}{2}} \left[I_{2n}, 0_{2n \times nK}\right], \quad L_b= Q_k^{\frac{1}{2}} \left[I_{2n}, 0_{2n \times (n+nK)}\right].
%\end{gather}
\begin{assumption}\label{ass:hurwitz}
The matrix $M_1$ is Hurwitz.
\end{assumption}
% \begin{assumption}\label{ass:observability} The pair $(L_{0,H}, M_3)$ is observable.
% \seb{what do we mean by observable here?}
% \dena{This was left from the previous versions, not needed here since we are considering finite horizon problems.}
% \end{assumption}
%\begin{assumption}\label{ass:detectability}
%The pair $(L_a, \mb{A}_0)$ is detectable, and the pair $(L_b, \mb{A}_k),\, k \in \mfK$, is detectable. 
%\end{assumption}
The following theorem links the infinite population equilibria to the finite
population case. %The analysis above\seb{but didn't we use this result to obtain some of the above results... something seems circular} leads to the following theorem, which we prove using convex analysis and an asymptotic MFG equilibrium analysis
%\seb{but didn't we use this result to obtain some of the above results... something seems circular} \dena{I revised the wording.}
\begin{theorem} \label{thm}
 Subject to \textit{Assumptions \ref{ass:IntialStateAss}-\ref{ass:hurwitz}},   the mean field equations \eqref{GSCMFGlatentCond1}-\eqref{GSCMFGlatentCond2} together with the system equations \eqref{GSCMFGMajorIncomDyn}-\eqref{GSCMFGminorComDynNonEst} and \eqref{GSCMFGlatent-majorcost}-\eqref{GSCMFGlatent-minorcost},  generate an infinite family of stochastic control laws  $\mathcal{U}_{MF}^{\infty, *}$, with finite sub-families $\mathcal{U}_{MF}^{N,*} ~\triangleq~ \{u^{i,*}_t;~  i \in \mathfrak{N} \}$, $1\leq N<\infty$,  given by \eqref{GSCMFGLatentSLQconvControl-major}-\eqref{GSCMFGLatentLQGoffset-major} and \eqref{GSCMFGLatentSLQconvControl}-\eqref{GSCMFGLatentLQGoffset}, such that
 \begin{itemize}
 
 \item[(i)] $ {\mathcal{U}}_{MF}^{\infty, *}$ yields a unique Nash equilibrium within the set of linear control laws $\mathcal{U}^{\infty}_{L}$ such that 
 \begin{align}
 J_i^{\infty}(u^{i,*}, u^{-i,*}) = \inf_{u^i \in \mathcal{U}_L^{\infty}} J_i^{\infty} (u^i, u^{-i,*}),
  \end{align}

\item[(ii)]  All agent systems $i \in \mathfrak{N}_0$, are second order stable in the sense that
\begin{equation}
 \sup_{t \in [0,T], i \in \mathfrak{N}_0} \mathbb{E} \Big \{( {\Vert x^i_t \Vert}^2 + {\Vert x^{(N)}_t \Vert}^2 +  {\Vert \bar{x}_t \Vert}^2  + {\Vert y_t \Vert}^2)  \Big \} < C,
\end{equation}
with $C$ independent of N.
 
 \item[(iii)] $\{ {\mathcal{U}}_{MF}^N;~ 1 \leq N < \infty \}$ yields a unique $\epsilon$-Nash equilibrium within the set of linear control laws $\mathcal{U}^N_{L}$ for all $\epsilon > 0$, i.e. for  all $\epsilon > 0$, there exists $N(\epsilon)$ such that for all $N \geq N(\epsilon)$
 \begin{equation}
 J_i^N(u^{i,*}, u^{-i,*}) - \epsilon \leq \inf_{u^i \in \mathcal{U}_L^N} J_i^N (u^i, u^{-i,*}) \leq J_i^N(u^{i,*}, u^{-i,*}),
  \end{equation}
  where  $J_i^N(u^{i,*}, u^{-i,*}) \rightarrow J_i^\infty(u^{i,*}, u^{-i,*})$,\, $i \in \mathfrak{N}_0$, as $N \rightarrow \infty$.
  \end{itemize} 
\end{theorem} 
\hfill $\square$
\begin{proof} First, we use the convex analysis method developed in \cite{FCJ-Convex2018} to obtain
the best response strategies (17)-(19) and (23)-(25); this proves part (i) of the theorem. Then, following the same lines as in the closed-loop and asymptotic equilibrium analysis of \cite[Sec. 5 \& 6]{Huang2010}, the agents are second order stable and the set of infinite population control actions yields an $\epsilon$-Nash equilibrium for the large population system which prove parts (ii) and (iii) of the theorem, respectively.
%\seb{needs rewording.}\dena{I revised this part.}

%For an overview of the convex analysis approach to optimization, see e.g., \cite{ConvecAnalysisBook1999} for static systems, and for the relationship between the G\^{a}teaux derivative of the cost functional of a system and its Hamiltonian see e.g., \cite{CarmonaBSDEsBook2016}. \cite{bank2017hedging} studies a stochastic tracking problem in finance using the convex analysis approach, while \cite{JaimungalPhil2018} investigates an algorithmic trading problem and obtains the optimal trading strategies for a large number of heterogeneous traders using the convex analysis approach. 

\underline{Part(i)}: Consider the major agent $\mc{A}_0$'s extended system. 
%to derive the major agent's optimal control action in the infinite population limit. 
Following the lines of the proof to Theorem 3.2 in \cite{FCJ-Convex2018},  the G\^{a}teaux derivative of the major agent's extended cost $\mc{D}J_0^{ex}(u^0)$ in the direction of $\omega^0_t \in \mc{U}^0$ is given by 
 \begin{multline}\label{majorGat}
\langle \mc{D}J_0^{ex}(u^0), \omega^0 \rangle = \mathbb{E} \bigg [\int_0^T \omega_t^{0\intercal} \bigg\{  \mb{N}_0^\intercal X^{0,u}_t + R_0 u^0_t  \\+ \mb{B}_0^\intercal\Big(e^{-\mb{A}_0^\intercal t} M_t^0 -\int_0^t e^{\mb{A}_0^\intercal(s-t) } \big( \mb{Q}_0 X^{0,u}_s + \mb{N}_0 u^0_s \big)ds \Big )  \bigg  \} dt\bigg ],
 \end{multline}
where the martingale $(M^0_t)_{t\in[0,T]}$ is given by  
\begin{align} \label{majorMartingale}
M_t^0 =  \mathbb{E} \Big[e^{\mb{A}_0^\intercal T} \mb{G}_0 X^{0,u}_T + \int_0^T e^{\mb{A}_0^\intercal s } ( \mb{Q}_0 X^{0,u}_s + \mb{N}_0 u^0_s )ds 
\, \Big | \, \mathcal{F}^{0,r}_t \Big].
 \end{align}
Given that \textit{Assumption \ref{ConvexityCondMajor}} holds, following the lines of the proof to Theorem 3.3 in \cite{FCJ-Convex2018},  the major agent $\mc{A}_0$'s optimal control action $u^{0,*}_t$  in the infinite population limit is
 \begin{equation} \label{majorCtrlRawForm}
 u^{0,*}_t = - R_0^{-1} \bigg [ \mb{N}_0^\intercal X^{0,*}_t
 +  \mb{B}_0^\intercal  \Big ( e^{-\mb{A}_0^\intercal t}M_t^0 -\int_0^t\!\! e^{\mb{A}_0^\intercal(s-t) }  \big(\mb{Q}_0 X^{0,*}_s + \mb{N}_0 u^{0,*}_s \big)ds\Big ) \bigg ],
 \end{equation}
which is obtained by setting \eqref{majorGat} to zero for all possible paths of  $\omega^0_t \in \mc{U}^0$.  
 
 Next,  define
 \begin{equation} \label{GSCMFGLatentSLQconvAdjoint-major}
 p^0_t := e^{-\mb{A}_0^\intercal t}M_t^0 - \int_0^t e^{\mb{A}_0^\intercal(s-t)} \big( \mb{Q}_0 X^{0,*}_s + \mb{N}_0 u^{0,*}_s  \big)ds,
 \end{equation}
 which is the adjoint process for the major agent's system in the stochastic maximum principle framework. Next, we adopt the ansatz
 \begin{equation} \label{GSCMFGLatentSLQansatz-major}
 p_t^0 = \Pi_0(t) X_t^{0,*} + s^0_t,
 \end{equation}
 whose substitution in  \eqref{majorCtrlRawForm} yields a linear state feedback form for the major agent's optimal control action, i.e. 
 \begin{equation}\label{GSCMFGLatentSLQconvControl-major3}
 u^{0,*}_t = - R_0^{-1} \big [ \mathbb{N}_0^\intercal X^{0,*}_t +  \mathbb{B}_0^\intercal \big(\Pi_0(t) X_t^{0,*} + s^0_t) \big) \big]. 
 \end{equation} 
 To determine $\Pi_0(t) \in \mdsR^{(2+K)n \times (2+K)n}$ and $s^0_t \in \mdsR^{(2+K)n}$, first both sides of \eqref{GSCMFGLatentSLQansatz-major} are differentiated and then \eqref{GSCMFGlatentMajorExtDyn} and \eqref{GSCMFGLatentSLQconvControl-major3} are substituted, which gives 
 \begin{multline}
 d{p}_t^0 = \Big[\big (\dot{\Pi}_0 + \Pi_0 \mathbb{A}_0  - \Pi_0 \mathbb{B}_0 R_0^{-1} \mathbb{N}_0^\intercal - \Pi_0 \mathbb{B}_0 R_0^{-1} \mathbb{B}_0^\intercal \Pi_0 \big) X^0_t dt \\ +\big(-\Pi_0 \mathbb{B}_0 R_0^{-1} \mathbb{B}_0^\intercal s^0_t  + \Pi_0 \mathbb{M}_t^0 \big)dt + d{s}^0_t\Big] + \Pi_0 {\Sigma}_0(t) dW^0_t.  \label{GSCMFGLatentSLQconvdiff1-major}
 \end{multline}
 Next, both sides of \eqref{GSCMFGLatentSLQconvAdjoint-major} are differentiated to yield
 \begin{equation} \label{GSCMFGLatentSLQconvAdjointDiff-major}
 dp_t^0 = (-\mb{A}_0^\intercal p^0_t - \mb{Q}_0X^0_t -\mb{N}_0 u^0_t)dt + e^{-\mb{A}_0^\intercal t}dM^0_t.
 \end{equation}
 The martingale representation theorem implies that $M^0_t$ can be written as  
 \begin{equation}\label{GSCMFGLatentMrtnglRep-major}
 M^0_t = M_0^0 + \int_0^t Z_s^0 dW_s^0, 
 \end{equation}
 where $Z_t^0$ is an $\mathcal{F}^{0,r}$-adapted process. 
%  Differentiating both sides of \eqref{GSCMFGLatentMrtnglRep-major} yields 
%   \begin{equation} \label{GSCMFGLatentMrtnglRepDiff-major}
%  dM_t^0 = Z_t^0 dW^0_t.
%   \end{equation}
Then, \eqref{GSCMFGLatentSLQconvControl-major3} and \eqref{GSCMFGLatentMrtnglRep-major} are substituted in \eqref{GSCMFGLatentSLQconvAdjointDiff-major} which implies
\begin{multline} \label{GSCMFGLatentSLQconvdiff2-major}
d{p}^0_t = \Big[(- \mb{Q}_0  + \mb{N}_0 R_0^{-1} \mb{N}_0^\intercal + \mb{N}_0 R_0^{-1} \mb{B}_0^\intercal \Pi_0 -\mb{A}_0^\intercal \Pi_0)X^{0,*}_t \\
+ (\mb{N}_0 R_0^{-1} \mb{B}_0^\intercal - \mb{A}_0^\intercal) s^0_t \Big]dt + q^0_t dW^0_t,
\end{multline} 
where $q_t^0 = e^{-\mb{A}_0^\intercal t} Z_t^0$.

Finally, \eqref{GSCMFGLatentSLQconvdiff1-major} and \eqref{GSCMFGLatentSLQconvdiff2-major} are equated which results in the deterministic Riccati equation
\begin{align}
\dot{\Pi}_0 + \Pi_0 \mb{A}_0 + \mb{A}_0^\intercal \Pi_0 - (\mb{B}_0^\intercal \Pi_0 + \mb{N}_0^\intercal)^\intercal R_0^{-1}( \mb{B}_0^\intercal \Pi_0 + \mb{N}_0^\intercal ) + \mb{Q}_0  &= 0, \\
\Pi_0(T) &= \mb{G}_0,
\end{align}
and the stochastic offset equation 
\begin{align}
%\begin{split}
 d{s}^0_t + \Big ( \big[(\mb{A}_0- \mb{B}_0R_0^{-1} \mb{N}_0^\intercal)^\intercal- \Pi_0 \mb{B}_0 R_0^{-1} \mb{B}_0^\intercal \big] s^0_t +  \Pi_0 \mb{M}^0_t \Big ) dt %\qquad %& \\
 + (\Pi_0 \Sigma_0 - q^0_t)dW^0_t =0, 
%\end{split}
 \\
s^0_T &= 0. 
\end{align}   
To derive the optimal control action for minor agent $\mc{A}_i$, $i \in \mathfrak{N}$,  as well as the corresponding Riccati and offset equations, a similar approach is followed. Following the lines of the proof to Theorem 3.2 in \cite{FCJ-Convex2018}, the G\^ateaux derivative of the extended cost functional $\mc{D}J^{ex}_k(u^{i}),\, k \in \mfK$, for minor agent $\mc{A}_i, i \in \mathfrak{N}$, is given by
\begin{multline}\label{minorGat}
\langle \mc{D}J_k^{ex}(u^i), \omega^i \rangle = \mathbb{E} \bigg [\int_0^T (\omega_t^i)^\intercal \bigg\{ \mb{N}_k^\intercal X^{i,u}_t + R_k u^i_t \\+ \mb{B}_k^\intercal\Big(e^{-\mb{A}_k^\intercal t} M_t^i -\int_0^t e^{\mb{A}_k^\intercal(s-t) } ( \mb{Q}_k X^{i,u}_s + \mb{N}_k u^i_s )ds \Big )  \bigg  \} dt\bigg ].
\end{multline}
where the martingale $(M^i_t)_{t\in[0,T]}$ is given by  
 \begin{equation}\label{martingaleMinor}
M_t^i :=  \mathbb{E} \Big[e^{\mb{A}_k^\intercal T} \mb{G}_k X^{i,u}_T +\int_0^T e^{\mb{A}_k^\intercal s } ( \mb{Q}_k X^{i,u}_s + \mb{N}_k u^i_s )ds \Big)  \,\Big |
\, \mathcal{F}^{i,r}_t \Big].
 \end{equation}

Given \textit{Assumption \ref{convexityCondMinor}}, following the lines of the proof to Theorem 3.3 in \cite{FCJ-Convex2018}, 
%\seb{you are referring to the theorem you are proving here... something seems wrong.}\dena{I corrected the wording.}  
the optimal control action $u^{i,*}_t$ for minor agent $\mc{A}_i$ in the infinite population limit is 
 \begin{equation} \label{minorControlRaw}
 u^{i,*}_t = - R_k^{-1} \Big [ \mb{N}_k^\intercal X^{i,*}_t 
  +  \mb{B}_k^\intercal  \Big (e^{-\mb{A}_k^\intercal t}M_t^i -\int_0^t e^{\mb{A}_k^\intercal(s-t) }  (\mb{Q}_k X^{i,*}_s + \mb{N}_k u^{i,*}_s )ds\Big ) \Big ],
 \end{equation}
 which is obtained by setting \eqref{minorGat} to zero for all possible paths of $\omega^i_t \in \mc{U}^i$. 
  
Next, define $p^i_t$ as 
 \begin{equation} \label{GSCMFGLatentSLQconvAdjoint}
 p^i_t = e^{-\mb{A}_k^\intercal t}M_t^i - \int_0^t e^{\mb{A}_0^\intercal (s-t)} ( \mb{Q}_k X^{i,*}_s + \mb{N}_k u^{i,*}_s  )ds,
 \end{equation}
 which corresponds to the adjoint process for the minor agent $\mc{A}_i$'s system in the stochastic maximum principle framework. Then adopting the ansatz 
 \begin{equation} \label{GSCMFGLatentSLQansatz}
 p_t^i = \Pi_k(t) X_t^{i,*} + s^{i,k}_t,
 \end{equation}
 and  substituting into \eqref{minorControlRaw} results in a linear state feedback form for $u^{i,*}_t$ as
 \begin{equation}\label{GSCMFGLatentSLQconvControlderiv}
 u^{i,*}_t = - R_k^{-1} \big [ \mathbb{N}_k^\intercal  X^{i,*}_t +  \mathbb{B}_k^\intercal  \big(\Pi_k(t) X_t^{i,*} + s^{i,k}_t \big) \big]. 
 \end{equation} 
To obtain $\Pi_k(t) \in \mdsR^{(3+K)n \times (3+K)n}$ and $s^{i,k}_t \in \mdsR^{(3+K)n}$, first both sides of \eqref{GSCMFGLatentSLQansatz} are differentiated and then \eqref{GSCMFGlatant-minorExtDyn} and \eqref{GSCMFGLatentSLQconvControl} are substituted which yields
 \begin{multline}
 d{p}_t^i = \Big[\big (\dot{\Pi}_k + \Pi_k \mathbb{A}_k  - \Pi_k \mathbb{B}_k R_k^{-1} \mathbb{N}_k^\intercal  - \Pi_k \mathbb{B}_k R_k^{-1} \mathbb{B}_k^\intercal  \Pi_k \big) X^{i,*}_t
 \\
 - \Pi_k\mathbb{B}_k R_k^{-1} \mathbb{B}_k^\intercal  s^{i,k}_t
 + \Pi_k \mathbb{M}_t^k + d{s}^{i,k}_t\Big] dt + \Pi_k {\Sigma}_k(t) dW^i_t.  \label{GSCMFGLatentSLQconvdiff1}
 \end{multline}
Next, both sides of \eqref{GSCMFGLatentSLQconvAdjoint} are differentiated yielding
 \begin{equation} \label{GSCMFGLatentSLQconvAdjointDiff}
 dp_t^i = (-\mb{A}_k^\intercal p^i_t - \mb{Q}_kX^{i,*}_t -\mb{N}_k^\intercal u^{i,*}_t)dt + e^{-\mb{A}_k^\intercal t}dM^i_t.
 \end{equation}
 According to the martingale representation theorem, the martingale $M^i_t$ admits the representation
 \begin{equation}\label{GSCMFGLatentMrtnglRep}
 M^i_t = M_0^i + \int_0^t Z_s^i dW_s^i, 
 \end{equation}
% or equivalently, when both sides of \eqref{GSCMFGLatentMrtnglRep} are differentiated, as 
% \begin{align} \label{GSCMFGLatentMrtnglRepDiff}
%  dM_t^i = Z_t^i dW^i_t,
% \end{align}
where $Z_t^i$ is an $\mc{F}^{i,r}$-adapted process. 

Then, \eqref{GSCMFGLatentSLQconvControl} and \eqref{GSCMFGLatentMrtnglRep} are substituted in \eqref{GSCMFGLatentSLQconvAdjointDiff} which gives
\begin{multline} \label{GSCMFGLatentSLQconvdiff2}
d{p}^i_t = \Big[(- \mb{Q}_k  + \mb{N}_k R_k^{-1} \mb{N}_k^\intercal + \mb{N}_k R_k^{-1} \mb{B}_k^\intercal \Pi_k -\mb{A}_k^\intercal \Pi_k)X^{i,*}_t \\+ (\mb{N}_k R_k^{-1} \mb{B}_k^\intercal - \mb{A}_k^\intercal) s^{i,k}_t \Big]dt + q^i_t dW^i_t,
\end{multline} 
where $q_t^i = e^{-\mb{A}_k^\intercal t} Z_t^i$.
Finally,  \eqref{GSCMFGLatentSLQconvdiff1} is equated with \eqref{GSCMFGLatentSLQconvdiff2} which yields 
\begin{align}
\dot{\Pi}_k + \Pi_k \mb{A}_k + \mb{A}_k^\intercal \Pi_k - (\mb{B}_k^\intercal \Pi_k + \mb{N}_k^\intercal)^\intercal R_k^{-1}( \mb{B}_k^\intercal \Pi_k + \mb{N}_k^\intercal ) + \mb{Q}_k  &= 0, \\
\Pi_k(T) &= \mb{G}_k,
\end{align} 
 \begin{align} 
 %\begin{split}
 ds^{i,k}_t + \Big ( \big[(\mb{A}_k- \mb{B}_kR_k^{-1} \mb{N}_k^\intercal)^\intercal- \Pi_k \mb{B}_k R_k^{-1} \mb{B}_k^\intercal \big] s^{i,k}_t + \Pi_k \mb{M}^k_t\Big )dt %\qquad %&
% \\
 + (\Pi_k \Sigma_k - q^i_t)dW^i_t =0, 
 %\end{split}
 \\
s^{i,k}_T &= 0,
\end{align}
$i \in \mathfrak{N}, k \in \mfK$.

\underline{Part (ii) \& (iii)}: Following the closed-loop and asymptotic equilibrium analysis of \cite[Sec. 5 \& 6]{Huang2010}, the set of control actions $\mathcal{U}_{MF}^{N,*} ~\triangleq~ \{u^{i,*}_t;~   i \in \mathfrak{N} \}$, $1 \leq N < \infty $, yields second order stability and an $\epsilon$-Nash equilibrium for the large population system given by \eqref{GSCMFGMajorIncomDyn}-\eqref{GSCMFGminorComDynNonEst} and \eqref{GSCMFGlatent-majorcost}-\eqref{GSCMFGlatent-minorcost}.
%\hfill $\square$
\end{proof}
%\section{Financial Example}

\section{Conclusions}\label{sec:conclusion}

In this paper, we introduced and formulated a new class of major-minor MFG systems motivated from financial and economic systems, such as those appearing in algorithmic trading and energy makets. In this novel setup, the major agent and each of the mass of minor agents interact with a common process, and this process also affects their cost functionals. The common process is influenced by (i) a latent process which is not observed, (ii) a common Wiener process, (iii) the major agent's state and control action, and (iv)  the average state and control action of all minor agents. Then, we developed a convex analysis method to establish the best trading strategies for all agents  which yield an $\epsilon$-Nash equilibrium. Our framework can be easily extended to the case where each agent's dynamics are also influenced by the common process.

 \bibliography{bib_Dena_19Jul19}

\begin{thebibliography}{10}

\bibitem{CainesEncyc2014}
P.~E. Caines, ``Mean field games,'' {\em Encyclopedia of Systems and Control},
  2019.

\bibitem{HuangCIS2006}
M.~Huang, R.~P. Malham\'e, and P.~E. Caines, ``Large population stochastic
  dynamic games: closed-loop {McKean-Vlasov} systems and the {Nash} certainty
  equivalence principle,'' {\em Communications in Information and Systems},
  vol.~6, no.~3, pp.~221--252, 2006.

\bibitem{HuangCDC2006}
M.~Huang, R.~P. Malham\'e, and P.~E. Caines, ``{N}ash certainty equivalence in
  large population stochastic dynamic games: Connections with the physics of
  interacting particle systems,'' in {\em Proceedings of the 45th {IEEE}
  Conference on Decision and Control ({CDC})}, (San Diego, {CA}),
  pp.~4921--4926, Dec. 2006.

\bibitem{HuangCDC2003a}
M.~Huang, P.~E. Caines, and R.~P. Malham\'e, ``Individual and mass behavior in
  large population stochastic wireless power control problems: centralized and
  {Nash} equilibrium solutions,'' in {\em Proceedings of the 42nd {IEEE}
  Conference on Decision and Control ({CDC})}, (Maui, {HI}), pp.~98--103, Dec.
  2003.

\bibitem{HuangCDC2003b}
M.~Huang, R.~P. Malham\'e, and P.~E. Caines, ``Stochastic power control in
  wireless communication systems: analysis, approximate control algorithms and
  state aggregation,'' in {\em Proceedings of the 42nd {IEEE} Conference on
  Decision and Control ({CDC})}, (Maui, {HI}), pp.~4231--4236, Dec. 2003.

\bibitem{HuangTAC2007}
M.~Huang, P.~E. Caines, and R.~P. Malham\'e, ``Large-population cost-coupled
  {LQG} problems with nonuniform agents: individual-mass behavior and
  decentralized $\epsilon$-{N}ash equilibria,'' {\em IEEE Transaction on
  Automatic Control}, vol.~52, no.~9, pp.~1560--1571, 2007.

\bibitem{CHM17}
P.~E. Caines, M.~Huang, and R.~P. Malham{\'e}, ``Mean field games,'' in {\em
  Handbook of Dynamic Game Theory} (T.~Ba\c{s}ar and G.~Zaccour, eds.),
  pp.~1--28, Berlin: Springer, 2017.

\bibitem{Lasry2006a}
J.~M. Lasry and P.~L. Lions, ``Jeux \`a champ moyen. i - le cas stationnaire,''
  {\em Comptes Rendus de l'Acad\'emie des Sciences}, vol.~343, pp.~619--625,
  2006.

\bibitem{Lasry2006b}
J.~M. Lasry and P.~L. Lions, ``Jeux \`a champ moyen. ii - horizon fini et
  contr\^ole optimal,'' {\em Comptes Rendus de l'Acad\'emie des Sciences},
  vol.~343, pp.~679--684, 2006.

\bibitem{Lasry2007}
J.~M. Lasry and P.~L. Lions, ``Mean field games,'' {\em Japanese Journal of
  Mathematics}, vol.~2, no.~1, pp.~229--260, 2007.

\bibitem{Huang2010}
M.~Huang, ``Large-population {LQG} games involving a major player: The {N}ash
  certainty equivalence principle,'' {\em SIAM Journal on Control and
  Optimization}, vol.~48, no.~5, pp.~3318--3353, 2010.

\bibitem{Huang2012}
S.~L. Nguyen and M.~Huang, ``Linear-quadratic-gaussian mixed games with
  continuum-parametrized minor players,'' {\em SIAM Journal on Control and
  Optimization}, vol.~50, no.~5, pp.~2907--2937, 2012.

\bibitem{NourianSiam2013}
M.~Nourian and P.~E. Caines, ``$\epsilon$-{N}ash mean field game theory for
  nonlinear stochastic dynamical systems with major and minor agents,'' {\em
  SIAM Journal on Control and Optimization}, vol.~51, no.~4, pp.~3302--3331,
  2013.

\bibitem{Kizilkale2013}
P.~E. Caines and A.~C. Kizilkale, ``Recursive estimation of common partially
  observed disturbances in {MFG} systems with application to large scale power
  markets,'' in {\em Proceedings of the 52nd IEEE Conference on Decision and
  Control (CDC)}, (Florence, Italy), pp.~2505--2512, Dec. 2013.

\bibitem{Kizilkale2014}
P.~E. Caines and A.~C. Kizilkale, ``Mean field estimation for partially
  observed {LQG} systems with major and minor agents,'' in {\em Proceedings of
  the 19th World Congress of the International Federation of Automatic Control
  (IFAC)}, (Cape Town, South Africa), pp.~8705--8709, Aug. 2014.

\bibitem{KizilkaleTAC2016}
P.~E. Caines and A.~C. Kizilkale, ``$\epsilon$-{N}ash equilibria for partially
  observed {LQG} mean field games with major player,'' {\em IEEE Transaction on
  Automatic Control}, vol.~62, no.~7, pp.~3225--3234, 2017.

\bibitem{SenCDC2014}
N.~\c{S}en and P.~E. Caines, ``Mean field games with partially observed major
  player and stochastic mean field,'' in {\em Proceedings of the 53rd IEEE
  Conference on Decision and Control (CDC)}, (Los Angeles, {CA}),
  pp.~2709--2715, Dec. 2014.

\bibitem{SenACC2015}
N.~\c{S}en and P.~E. Caines, ``$\epsilon$-{N}ash equilibria for a partially
  observed mean field game with major player,'' in {\em Proceedings of the 2015
  American Control Conference (ACC)}, (Chicago, IL), pp.~4791--4797, July 2015.

\bibitem{SenSiam2016}
N.~\c{S}en and P.~E. Caines, ``Mean field game theory with a partially observed
  major agent,'' {\em SIAM Journal on Control and Optimization}, vol.~54,
  no.~6, pp.~3174--3224, 2016.

\bibitem{FirooziCDC2015}
D.~Firoozi and P.~E. Caines, ``$\epsilon$-{N}ash equilibria for partially
  observed {LQG} mean field games with major agent: Partial observations by all
  agents,'' in {\em Proceedings of the 54th IEEE Conference on Decision and
  Control (CDC)}, (Osaka, {Japan}), pp.~4430--4437, Dec. 2015.

\bibitem{JaimungalBook2015}
{\'A. Cartea and S. Jaimungal and J. Penalva}, {\em Algorithmic and
  high-frequency trading}.
\newblock {United Kingdom }: {Cambridge University Press}, 2015.

\bibitem{Mojtaba2015}
X.~Huang, S.~Jaimungal, and M.~Nourian, ``Mean-field game strategies for
  optimal execution,'' {\em Applied Mathematical Finance}, vol.~26, no.~2,
  pp.~153--185, 2019.

\bibitem{FirooziCDC2016}
D.~Firoozi and P.~E. Caines, ``Mean field game $\epsilon$-{N}ash equilibria for
  partially observed optimal execution problems in finance,'' in {\em
  Proceedings of the 55th IEEE Conference on Decision and Control (CDC)}, (Las
  Vegas, {NV}, {USA}), pp.~268--275, Dec. 2016.

\bibitem{FirooziCDC2017}
D.~Firoozi and P.~E. Caines, ``An optimal execution problem in finance
  targeting the market trading speed: an {MFG} formulation,'' in {\em
  Proceedings of the 56th IEEE Conference on Decision and Control (CDC)},
  (Melbourne, Australia), pp.~7--14, Dec. 2017.

\bibitem{FirooziIFAC2017}
D.~Firoozi and P.~E. Caines, ``An optimal execution problem in finance with
  acquisition and liquidation objectives: an {MFG} formulation,'' {\em
  IFAC-PapersOnLine}, vol.~50, no.~1, pp.~4960 -- 4967, 2017.

\bibitem{FirooziISDG2017}
D.~Firoozi and P.~E. Caines, ``The execution problem in finance with major and
  minor traders: A mean field game formulation,'' in {\em Annals of the
  International Society of Dynamic Games (ISDG): Advances in Dynamic and Mean
  Field Games}, vol.~15, pp.~107--130, Birkh\"auser Basel, 2017.

\bibitem{FirooziPakniyatCainesCDC2017}
D.~Firoozi, A.~Pakniyat, and P.~E. Caines, ``A mean field game - hybrid systems
  approach to optimal execution problems in finance with stopping times,'' in
  {\em Proceedings of the 56th IEEE Conference on Decision and Control (CDC)},
  (Melbourne, Australia), pp.~3144--3151, Dec. 2017.

\bibitem{Cardaliaguet2018}
P.~Cardaliaguet and C.-A. Lehalle, ``Mean field game of controls and an
  application to trade crowding,'' {\em Mathematics and Financial Economics},
  vol.~12, no.~3, pp.~335--363, 2018.

\bibitem{Lehalle2019}
C.-A. Lehalle and C.~Mouzouni, ``A mean field game of portfolio trading and its
  consequences on perceived correlations,'' {\em arXiv preprint
  arXiv:1902.09606}, 2019.

\bibitem{JaimungalPhil2018}
P.~Casgrain and S.~Jaimungal, ``Algorithmic trading with partial information: A
  mean field game approach,'' {\em arXiv}, 2018.

\bibitem{casgrain2018mean}
P.~Casgrain and S.~Jaimungal, ``Mean-field games with differing beliefs for
  algorithmic trading,'' {\em arXiv preprint arXiv:1810.06101}, 2018.

\bibitem{Horst2018}
G.~Fu, P.~Graewe, U.~Horst, and A.~Popier, ``A mean field game of optimal
  portfolio liquidation,'' {\em arXiv preprint arXiv:1804.04911}, 2018.

\bibitem{Aid2017}
R.~A\"{i}d, M.~Basei, and H.~Pham, ``The coordination of centralised and
  distributed generation,'' {\em arXiv preprint arXiv:1705.01302}, 2017.

\bibitem{JaimungalPhil2016}
P.~Casgrain and S.~Jaimungal, ``Trading algorithms with learning in latent
  alpha models,'' {\em SSRN}, 2016.

\bibitem{CarmonaDelarue2016}
R.~Carmona, F.~Delarue, and D.~Lacker, ``{Mean field games with common
  noise},'' {\em Annals of Probability}, vol.~44, pp.~3740--3803, 2016.

\bibitem{FJCCDC2018}
D.~Firoozi, S.~Jaimungal, and P.~E. Caines, ``Mean field game systems with
  common and latent processes,'' in {\em Proceedings of the 57th IEEE
  Conference on Decision and Control (CDC)}, ({Miami Beach}, {FL}, {USA}),
  pp.~5500--5505, Dec. 2018.

\bibitem{XYZBook1999}
{J. Yong and X. Y. Zhou}, {\em Stochastic controls: Hamiltonian systems and
  {HJB} Equations}.
\newblock { New York}: {Springer-Verlag }, 1999.

\bibitem{Wonham1964}
W.~M. Wonham, ``Some applications of stochastic differential equations to
  optimal nonlinear filtering,'' {\em Journal of the Society for Industrial and
  Applied Mathematics Series A Control}, vol.~2, no.~3, pp.~347--369, 1964.

\bibitem{CarmonaBSDEsBook2016}
R.~Carmona, {\em Lectures on BSDEs, stochastic control, and stochastic
  differential games with financial applications}.
\newblock Philadelphia, PA: Society for Industrial and Applied Mathematics,
  2016.

\bibitem{FCJ-Convex2018}
D.~Firoozi, S.~Jaimungal, and P.~E. Caines, ``Convex analysis for {LQG} systems
  with applications to major minor {LQG} mean field game systems,'' {\em
  arXiv}, 2018.

\end{thebibliography}
  \bibliographystyle{ieeetr}

\end{document}